\documentclass[a4paper,12pt]{article}

\newcommand{\kmcomment}[1]{}

\newcommand{\ds}{\ensuremath{\displaystyle }}

\renewcommand{\c}{\sigma}

\newcommand{\mR}{\ensuremath{\mathbb{R}}} %Use in Math-mode
\newcommand{\mZ}{\ensuremath{\mathbb{Z}}} %Use in Math-mode
\newcommand{\ANYb}[2]{{#1}_{#2}}

\newcommand{\yb}[1]{y_{#1}}
\newcommand{\zb}[1]{z_{#1}}
\newcommand{\wb}[1]{w_{#1}}

\newcommand{\frakg}{\mathfrak{g}}

\newcommand{\frakh}{\mathfrak{h}}
\newcommand{\frakhN}[1]{\mathfrak{h}_{#1}}

\newcommand{\tbdl}[1]{\mathrm{T}(#1)}

\newcommand{\cbdl}[1]{\mathrm{T}^{*}(#1)}

\newcommand{\Sbt}[2]{[#1,#2]}                
\newcommand{\parity}[1]{(-1)^{#1}}

\newcommand{\Akt}[2]{[#1,#2]}

\newcommand{\Bkt}[2]{\{#1,#2\}_{\pi}}

\newcommand{\mw}{\mywedge}
\newcommand{\we}{\wedge}

\renewcommand{\dim}{\textrm{dim}}

\newcommand{\mywedge}{\bigtriangleup} 

\renewcommand{\[}{$$} \renewcommand{\]}{$$}

\usepackage{ntheorem} 
{%
\theoremheaderfont{\bfseries}
\theorembodyfont{\rmfamily}

\newtheorem{prop}{Proposition}[section]

\newtheorem{remark}{Remark}[section]

\newtheorem{thm}{Theorem}[section]
\newtheorem*{thm-none}{Theorem}[section]

}

\renewcommand{\[}{$$} \renewcommand{\]}{$$}

\usepackage{fancybox,fancyvrb,shortvrb}
\usepackage{colortbl,hhline}
\usepackage{mathtools}
\usepackage{listings,jlisting}

\usepackage{setspace,Xentry}

\usepackage{graphicx}
\usepackage{hyperref}
\usepackage{maple}
\usepackage[utf8]{inputenc}
\usepackage[svgnames]{xcolor}
\usepackage{amsmath,amsfonts,newtxmath,newtxtext}
\usepackage{breqn}
\usepackage{textcomp}
\usepackage{pdfpages}
\usepackage{tikz} % after xcolor
\usetikzlibrary{positioning,calc,quotes,angles}
\usetikzlibrary{backgrounds}

\usepackage{pxpgfmark} %remember picture を可能にする
 \RequirePackage{listings} 
 \lstnewenvironment{kmverb}
    { \lstset{xleftmargin=10mm, breaklines=true}}
    { }

 \lstnewenvironment{kmverbb}[1]
    { \lstset{xleftmargin=#1, showspaces=false, breaklines=true}}
    { }
 \lstnewenvironment{kmverbN10}
    { \lstset{xleftmargin=-10mm, breaklines=true}}
    { } 
 \lstnewenvironment{kmverbZero}
    { \lstset{xleftmargin=0mm, breaklines=true}}
    { }

 \lstnewenvironment{kmverb*}
    { \lstset{xleftmargin=10mm, showspaces=true, breaklines=true}}
    { }

\usepackage{wrapfig}

\usepackage{subsupscripts}
\DefineShortVerb{\!} %\UndefineShortVerb{\!}

\numberwithin{equation}{section}
\allowdisplaybreaks
\usepackage{anysize}
\marginsize{10mm}{10mm}{5mm}{5mm}

\usepackage[enableskew,vcentermath]{youngtab} % \tiny \yng(4,1,1)

\usepackage{comment}
\usepackage{here}

\title{%An introduction of 
Application of superalgebra homology groups to distinguish Engel-like
structures}

\author{Kentaro Mikami\footnote{
Akita University, Japan (mikami@math.akita-u.ac.jp) \quad  
This work is supported by JSPS KAKENHI Grant Number 22K03306. 
} 
\and Tadayoshi Mizutani\footnote{Saitama University, Japan} \and
Hajime Sato \footnote{Nagoya University, Japan} 
}

\allowdisplaybreaks
\AtBeginDocument{\parindent=0pt \setlength{\headheight}{16pt}
  \abovedisplayskip     =0.5\abovedisplayskip
  \abovedisplayshortskip=0.5\abovedisplayshortskip
  \belowdisplayskip     =0.5\belowdisplayskip
  \belowdisplayshortskip=0.5\belowdisplayshortskip}

\begin{document}
\thispagestyle{plain}

\maketitle

\section*{Abstract}
For a long time, tangent bundle of a manifold and the direct sum of
multi-vector fields is a prototype of \(\mZ\)-graded superalgebra with the
Schouten bracket, and the grading is defined by 
\(a-1\) for \(a\)-multivector fields. 
Recently, we noticed the direct sum of differential forms
of a manifold has \(\mZ\)-graded superalgebra structure with a super bracket 
\( \Sbt{A}{B} = \parity{a} d ( A \wedge B )\) for \(a\)-form \(A\),  
and the grading of \(a\)-form is given by \(-a -1\).  By Lie derivation, we
also see that \(\sum\limits_{p=0} ^{n}\Lambda^{p} \cbdl{M} \oplus \tbdl{M}\)
has 
\(\mZ\)-graded Lie superalgebra structure, which is a ``super'' superalgebra of the two
superalgebra described above. 

For a given \(\mZ\)-graded Lie superalgebra, there is a notion of the
weighted (co)chain complex and (co)homology groups. In general, those
objects are infinite dimensional and hard to understand the entire
properties.  If we restrict the manifold above to a finite dimensional Lie
group, and multivector fields and differential forms to (left) invariant
fields and forms, then (co)chain spaces become finite dimensional.  Still
studying \(\mZ\)-graded Lie super algebra, which is a generalization of Lie
algebra, we encounter complicated manipulation of sign changes depending on
each degree. In this note, we introduce our trial of using Maple software to
gather possible equations of ``Engel like'' Lie
algebra structures by Jacobi identify, then simply solve them and get six
possibilities.   
Getting weighted chain complex and manipulate its homology groups,    
we claim that those six possibilities are not isomorphic as Lie algebras in
generic. 
\kmcomment{
At first, we ask Maple to gather possible equations of ``Engel like'' Lie
algebra structures on 4-dimensional space by Jacobi identify.   
Next, we ask Maple to solve the above possible equations.  As you see below, 
Maple shows six possibilities, which we call ourType[i], \(i=1,\ldots,6\).    
They  give  4-dimensional Lie algebras.  A simple question is those are
distinguished each other, namely not Lie algebra isomorphic or not.  
For this purpose, we check the weighted homology groups for three kinds of
superalgebras. This is our strategy in this note. 
Our claim is those six ourTypes are not isomorphic as Lie algebras in
generic. 
}

\section*{Keywords}

\(\mZ\)-graded Lie super algebra (which is a generalization of Lie algebra),  super bracket, the Schouten bracket, chain spaces, the boundary operator, 
the homology groups, the Betti numbers, and the Euler number.  

%\newpage
\newcommand{\CC}[1]{C_{#1}}

\section{Introduction}

For a 4-dimensional manifold \(M\), if rank 2 
distribution \(D\) satisfies \(D[2] := D + [D,D]\) has rank 3 and 
 \(D[3] := D[2] + [D[2],D[2]]\) has rank 4, then the pair \((M,D)\) is called an Engel
 structure.  Here, we apply the above notion for 4-dimensional Lie algebra
 \(\frakg\).  
The pair \( ( \frakg, D )\) is a Engel-type if \(D\) is a 2-dimensional
 subspace of \(\frakg\),  and 
 \(\dim D[2] = 3\) where \(D[2] := D + [D,D]\), and  
 \(\dim D[3] = 4\) where \(D[3] := D[2] + [D[2],D[2]]\). 
From the dimensional restriction, 
we see a specific basis \(\yb{1}, \yb{2}\) so that \( \Sbt{\yb{1}}{\yb{2}}=
 \yb{3} \in [D,D]\). We choose \(\yb{4}\) by   \( \Sbt{\yb{1}}{\yb{3}}=
 \yb{4}\).  Thus, 
\begin{alignat*}{5}
 \Sbt{\yb{1}}{\yb{2}} &=  \yb{3}\;,  
 \qquad &  
 \Sbt{\yb{2}}{\yb{1}} &= - \yb{3} 
%\\
& \qquad 
 \Sbt{\yb{1}}{\yb{3}} &= \yb{4}\;,  
 \quad &
 \Sbt{\yb{3}}{\yb{1}} &= -  \yb{4} 
\\
 \Sbt{\yb{1}}{\yb{4}} &= \sum_{i} C_{1,4,i} \yb{i}\;,  
 \quad &
 \Sbt{\yb{4}}{\yb{1}} &= - \sum_{i} C_{1,4,i} \yb{i} 
%\\
& \qquad 
 \Sbt{\yb{2}}{\yb{3}} &= \sum_{i} C_{2,3,i} \yb{i}\;,  
 \quad &
 \Sbt{\yb{3}}{\yb{2}} &= - \sum_{i} C_{2,3,i} \yb{i} 
\\
 \Sbt{\yb{2}}{\yb{4}} &= \sum_{i} C_{2,4,i} \yb{i}\;,  
 \quad &
 \Sbt{\yb{4}}{\yb{2}} &= - \sum_{i} C_{2,4,i} \yb{i} 
%\\
& \qquad 
 \Sbt{\yb{3}}{\yb{4}} &= \sum_{i} C_{3,4,i} \yb{i}\;,  
 \quad &
 \Sbt{\yb{4}}{\yb{3}} &= - \sum_{i} \CC{3,4,i} \yb{i} 
 \; . 
\end{alignat*}

We let know Maple2021 about the relations above and check Jacobi identities
by our maple script \texttt{Engel-try-1.mpl}. 
\begin{spacing}{0.92}
%\lstinputlisting{Engel-try-1-crone.mpl}
\begin{lstlisting}
(* Jun 19, 2022 dim gg = 4; D is a 2dim plane of gg.
D^{2} := D + [D, D]: D^{3} := D^{2} + [ D^{2}, D^{2} ]: 
If dim D^{2} = 3 and dim D^{3} = 4, ( gg , D ) is called an Engel structure.
Now want to find concrete Lie algebra having Engel structure and whose Engel's.  *)

C := table(): 
with(difforms): defform(y=1, C=0): 
n := 4: ybase := [seq( y[i],i=1..n) ]: gBase := ybase: 
read"Sbt-renew-v1.mpl": # SbtN and those are in BIG.mla   
# with(`km/LieAlg/Sbt-renew-v1.mpl`): 
kmread("km/LieAlg/Sbt-renew-v1.mpl"): 
#For test kmread("X/Y/zz.mpl"); 

for i to n do SbtT[ y[i], y[i] ] := 0 od: 
SbtT[y[1],y[2]] := y[3]: SbtT[y[2],y[1]] := - y[3]: 
SbtT[y[1],y[3]] :=  y[4] : SbtT[y[3],y[1]] := - y[4] : 

SbtT[y[1],y[4]] := add(C[1,4,j]*y[j] , j=1..n): 
SbtT[y[4],y[1]] := -add(C[1,4,j]*y[j] , j=1..n): 
myVar := NULL: myVar ,= seq(  C[1,4,j], j=1..n): 

for i from 2 to n do for k from i+1 to n do 
SbtT[y[i],y[k]] := add( C[i,k,j]*y[j] , j=1..n) : 
SbtT[y[k],y[i]] := -add( C[i,k,j]*y[j] , j=1..n) : 
myVar ,= seq(  C[i,k,j],  j=1..n): od od: 

myJac := proc( A,B,C ) local a,b,c;
SbtN(A, SbtN(B,C)) + SbtN(B, SbtN(C,A)) + SbtN(C, SbtN(A,B)) : end proc:  
mySkew  := proc( A,B ) local a,b,c; SbtN(A, B) + SbtN(B, A) : end proc:  
(* ukeComm := NULL:  
for i to n do for j to n do ukeComm ,= mySkew( y[i], y[j] ) od od: *)

ukeJ := NULL:  for i to n do for j from i+1 to n do for k from j+1 to n do 
ukeJ ,= myJac ( y[i], y[j], y[k] ) od od od: 
myEqn := NULL:
for i to n do myEqn ,= seq( diff( ukeJ[i], y[j] ), j=1..n ) od: 
mySol := solve({myEqn}, {myVar}); 
#NG mySolK := map( simplify, mySol ): 
mySolK := map( curry( map, simplify ), [ mySol ] ): 
mySolTools := SolveTools[PolynomialSystem]({myEqn}, {myVar}); 
mySolToolsK := map( curry(map, simplify), [mySolTools] ): 
for j to nops( mySolK ) do 
mySolInd,mySolDep := selectremove( i-> lhs(i)= rhs(i), mySolK[j] ):
ourSolInd[j] := map( lhs, mySolInd ):
ourSolDep[j] := mySolDep  :
mySolToolInd, mySolToolDep := selectremove( i-> lhs(i)= rhs(i), mySolToolsK[j] ):
ourSolToolInd[j] := map( lhs, mySolToolInd ):
ourSolToolDep[j] := mySolToolDep: 
od: 
(* UKE := table(sparse=NULL): SbtTT := table(): 
for i to nops([mySol]) do print(i,"", mySol[i]):
for j to n do for k from j+1 to n do 
    UKE[i] ,= SbtTT [y[j],y[k]] = simplify( SbtT[y[j],y[k]], mySol[i] ) od od: 
print(UKE[i]): od: *) 
# assign( UKE[1] ): BktT[ y[i], y[j] ]: 
# save mySol, SbtT, UKE, "Engel-try-1-out2.txt": 
# save ourSolToolInd, ourSolToolDep, ourSolInd, ourSolDep, SbtT, "Engel-try-1-out4.txt": 
(* Important caution. ourSol and ourSotTool are not exactly equal. 1 to 4 are equal but 
5 and 6 are swapped. So, we do use one way,  do not mixed using.  Jul 14, 2022 *)
(* i := 1: while( X[i] <> [] ) do some job : i := i+1 od:  *)

\end{lstlisting}
\end{spacing}
We show how Maple works well. 

%\newpage 

% \includegraphics{~/eng-try-1-out.pdf}
\kmcomment{
\includepdf[pages={1,2}, scale=1,offset = 0mm 0mm, trim = 0mm 0mm 0mm
0mm, frame=false]{eng-try-1-out.pdf}

\vspace*{-10cm}
\begin{minipage}{\textwidth}
\includepdf[pages={3}, scale=1, frame=false]{cropOut.pdf}
\end{minipage}
  }

\begin{minipage}{\textwidth}
%\includepdf[pages={1}, scale=0.9, offset = 0mm -115mm, 
\includepdf[pages={1}, scale=0.9, offset = 0mm -125mm, 
trim = 0mm 0mm 0mm
0mm, pagecommand={\thispagestyle{plain}}, frame=false]
{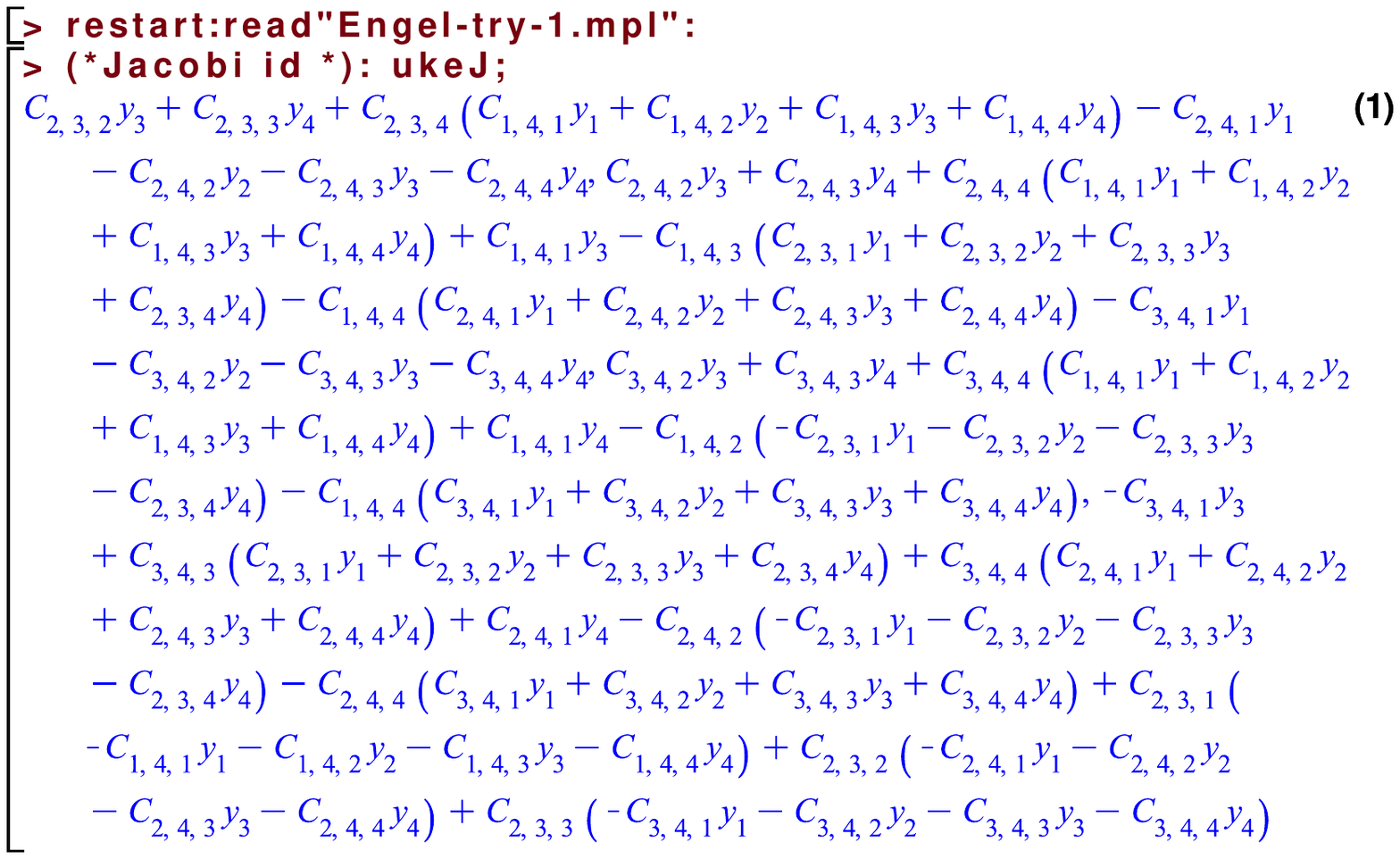}
%{/home/mikami/Engel-try-out1.pdf}
\end{minipage}

\newpage 

\begin{minipage}{\textwidth}
\includepdf[pages={1}, scale=0.9,offset = 0mm 15mm, trim = 0mm 0mm 0mm
0mm, pagecommand={\thispagestyle{plain}}, frame=false]
%{/home/mikami/Engel-try-out2.pdf}
{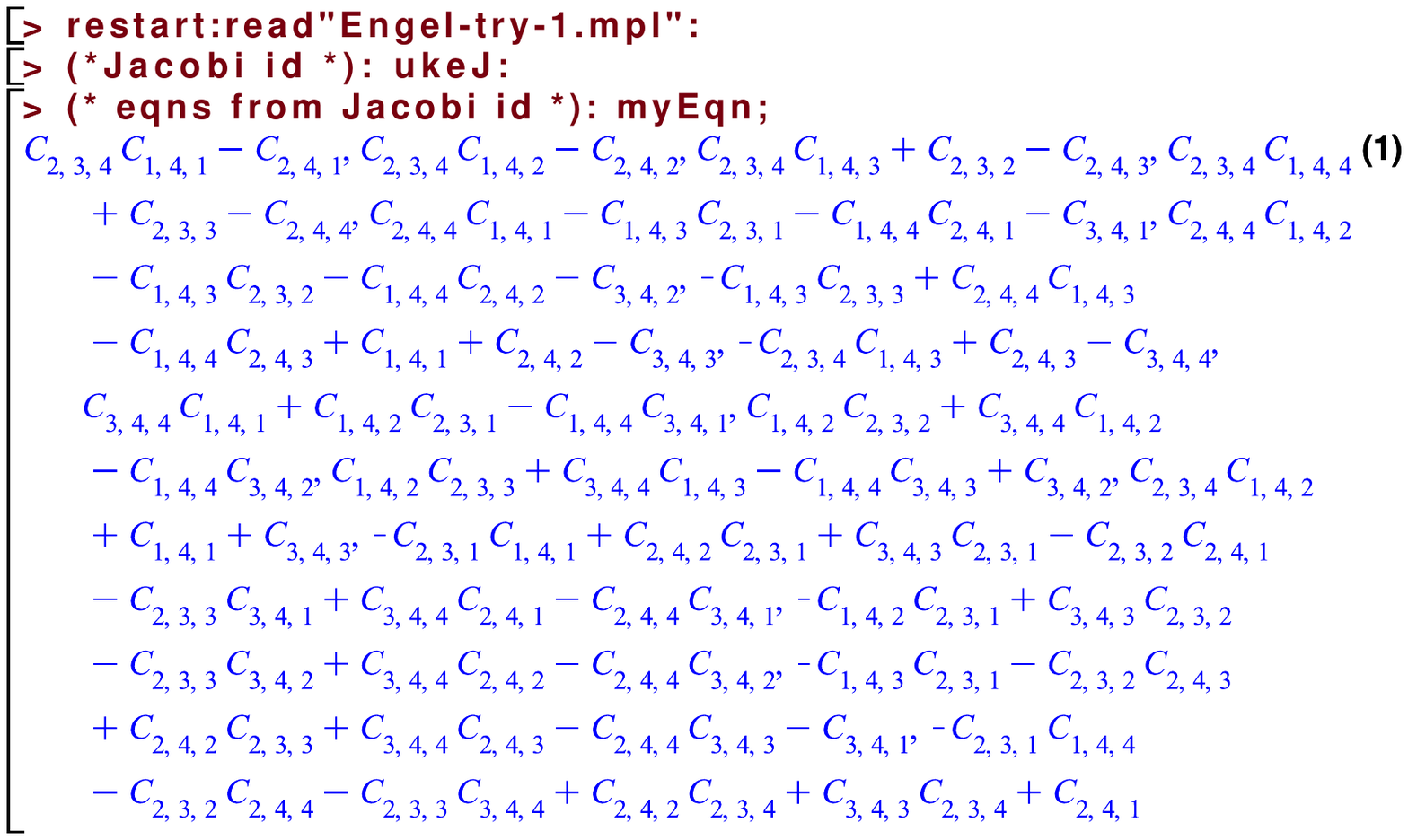}
\end{minipage}

\begin{minipage}{\textwidth}
\includepdf[pages={1}, scale=0.9,offset = 0mm -80mm, trim = 0mm 0mm 0mm
0mm, pagecommand={\thispagestyle{plain}}, frame=false]
%{/home/mikami/Engel-try-out3.pdf}
{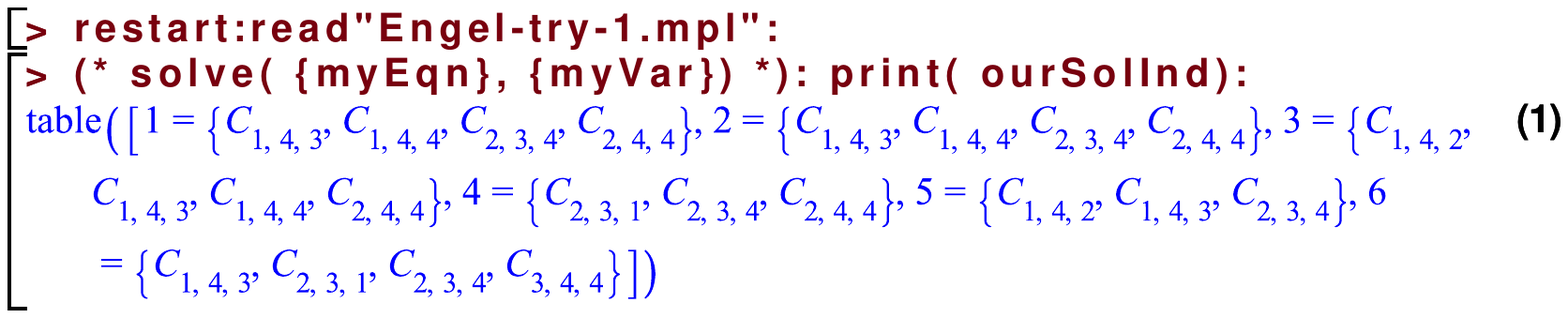}
\end{minipage}

\begin{minipage}{\textwidth}
\includepdf[pages={1}, scale=0.9,offset = 0mm -115mm, trim = 0mm 0mm 0mm
0mm, pagecommand={\thispagestyle{plain}}, frame=false]
%{/home/mikami/Engel-try-out4.pdf}
{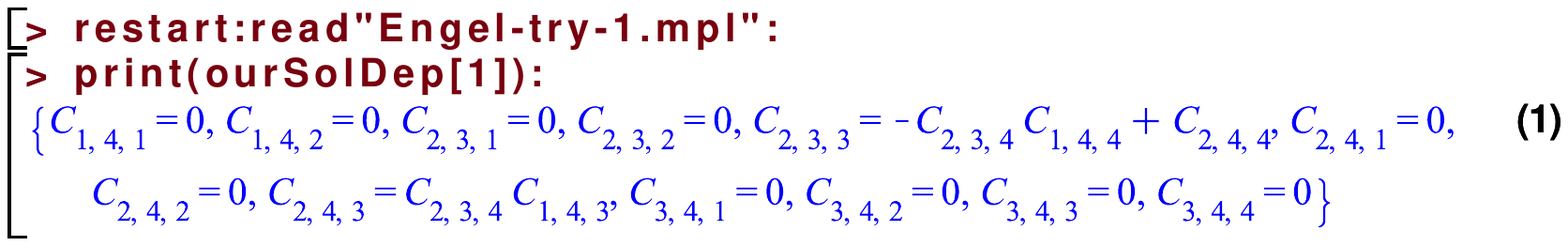}
\end{minipage}

\begin{minipage}{\textwidth}
\includepdf[pages={1}, scale=0.9,offset = 0mm -140mm, trim = 0mm 0mm 0mm
0mm, pagecommand={\thispagestyle{plain}}, frame=false]
%{/home/mikami/Engel-try-out6.pdf}
{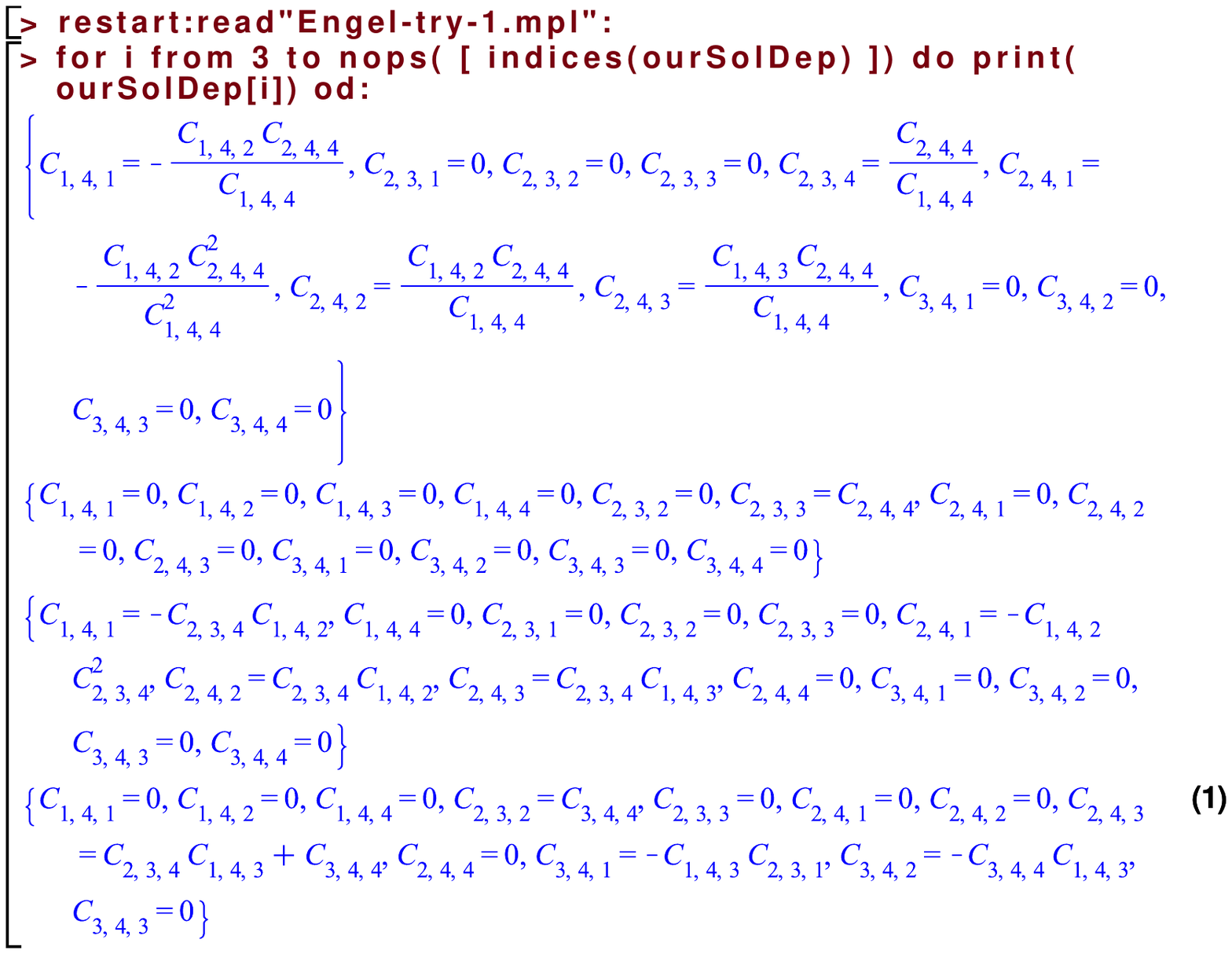}
\end{minipage}

\newpage

\begin{minipage}{\textwidth}
\includepdf[pages={1}, scale=0.9,offset = 0mm 5mm, trim = 0mm 0mm 0mm
0mm, pagecommand={\thispagestyle{plain}}, frame=false]
%{/home/mikami/Engel-try-out5.pdf}
{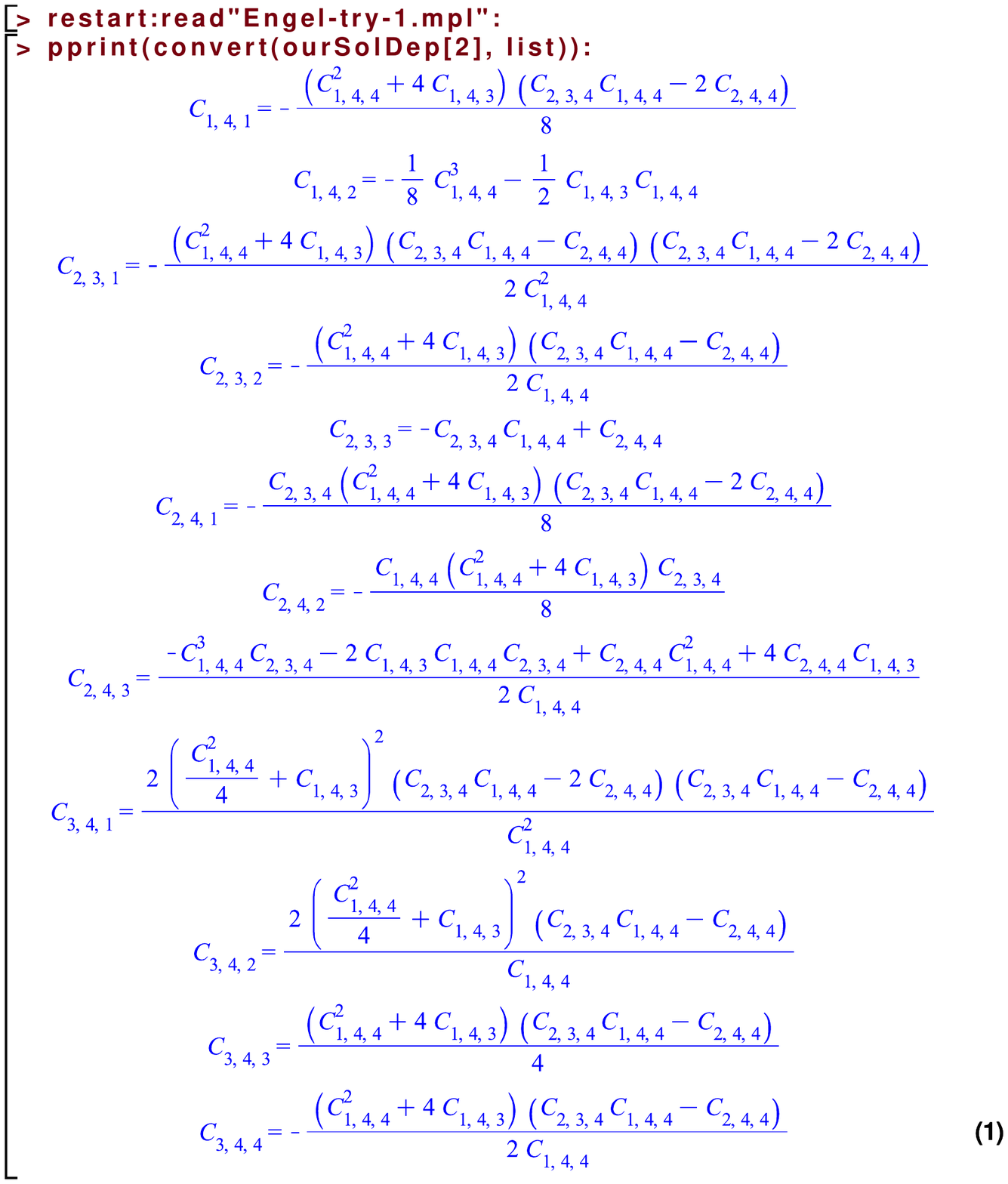}
\end{minipage}

\kmcomment{ 

\begin{minipage}{\textwidth}
\includepdf[pages={1}, scale=0.7,offset = 0mm -25mm, trim = 0mm 0mm 0mm
0mm, pagecommand={\thispagestyle{plain}}, frame=false]
{July17-Eng-try-1-ourSolDep.pdf}
\end{minipage}

\newpage 

\begin{minipage}{\textwidth}
\includepdf[nup=2x1,pages={2,3}, scale=1.4,offset = 40mm 40mm, trim = 0mm 0mm 0mm
0mm, pagecommand={\thispagestyle{plain}}, frame=false]
{July17-Eng-try-1-ourSolDep.pdf}
\end{minipage}

}

\kmcomment{  
%\begin{minipage}{\textwidth}
\includepdf[pages={3}, scale=1.0,trim = 0mm 250mm 0mm 10mm, frame=false,
clip,  offset=0 130mm, 
pagecommand={\thispagestyle{plain}} ]
{July17-Eng-try-1-ourSolDep.pdf}
%\end{minipage}
}

% \newpage 

\vspace*{180mm}

As you expect, in the original \texttt{Engel-try-1.mpl} if we declare
\texttt{save} at the almost end of the file, we get an output file
\texttt{Engel-try-1-out4.txt}.
Preparing the general form of Lie bracket candidate of Engel-type, we gather
Jacobi conditions \texttt{myEqn}. Solving them, we have 6 cases, and denote
them as \texttt{ourSolInd[i]} and \texttt{ourSolDep[i]} for \(i=1,
\ldots, 6\). As you see, the independent variable \(\CC{1,4,4}\) of the
second and third cases is assumed to be non-zero.   \texttt{ourSolDep[2]}
looks rather complicated comparing the others. Our main purpose of this note
is to distinguish those 6 cases are Lie algebra isomorphic or not by
watching their homology groups in superalgebra sense. 
Our strategy is simple, that is, if two Lie algebras \(\frakg\) and
\(\frakh\) are isomorphic as Lie algebras, their superalgebras are
isomorphic, and so weighted homology groups coincide  for each primary weight and the superalgebra type.  
Contrary, if the corresponding homology groups are not equal 
 for some same type of superalgebras and for some primary weight, 
 then the original Lie algebras are not isomorphic.

\section{Application of super homology of tangent type
 \(\ds \mathop{\oplus}_{\ell=1} ^{4} \Lambda ^{\ell} \frakg\) }
We have six types after solving the Jacobi identity conditions. 
In the type 2 and 3, some coefficients are fractions of denominator 
\(\CC{1, 4,4}\), this means it is non-zero. Looking at those six types, the
type 2 seems to be more complicated. 

We have a super homology group theory for superalgebra, especially derived
from Lie algebras which come from Engel like structure as we got above. 
We ask Maple to manipulate superalgebra homology groups with the primary
weight 0, 1, or 2. In the following, we show the outputs by Maple for the
six types in "generic". In the weight 0 case, the super homology group is
just the Lie algebra homology group and it is natural to treat 0-th chain
space \( \frakg^{0} = \mR\), and modify the output by Maple and we conclude
the Euler number is 0 for the weight 0 case, too.       
% So far, we are not able to manage vain spaces between pdf inputs. 

\begin{table}
\centering
%\begin{center}
\begin{tabular}{c@{\hspace{10mm}}c@{\hspace{10mm}}c}
{Weight 0, Common Headers} 
& 
{Weight 1, Common Headers} 
&
{Weight 0, Common Headers} 
\\
\( 
\left[
\begin{array}{*{6}{c}}
\text{m} & 0 & 1 & 2 & 3 & 4 \\
\text{SpaD} & 1 & 4 & 6 & 4 & 1 \\
\end{array}
\right]
\)
& 
\( 
\left[
\begin{array}{*{6}{c}}
\text{m} & 1 & 2 & 3 & 4 & 5 \\
\text{SpaD} & 6 & 24 & 36 & 24 & 6 \\
\end{array}
\right]
\)
& 
\(\left[ 
\begin{array}{*{7}{c}}
\text{m} & 1 & 2 & 3 & 4 & 5 & 6 \\
\text{SpaD} & 4 & 37 & 108 & 142 & 88 & 21 \\
\end{array}\right]
\)
\\
{Weight 0 of Engel ourType 1} 
& 
{Weight 1 of Engel ourType 1} 
&
{Weight 2 of Engel ourType 1} 
\\
\( 
\left[
\begin{array}{*{6}{c}}
\text{KerD} & 1 & 4 & 4 & 1 & 0 \\
\text{Bett} & 1 & 2 & 1 & 0 & 0 \\
\end{array}
\right]
\)
& 
\( 
\left[
\begin{array}{*{6}{c}}
\text{KerD} & 6 & 19 & 19 & 6 & 0 \\
\text{Bett} & 1 & 2 & 1 & 0 & 0 \\
\end{array}
\right]
\)
& 
\(\left[ 
\begin{array}{*{7}{c}}
\text{KerD} & 4 & 33 & 76 & 68 & 21 & 0 \\
\text{Bett} & 0 & 1 & 2 & 1 & 0 & 0 \\
\end{array}\right]
\)
\\
{Weight 0 of Engel ourType 2} 
& 
{Weight 1 of Engel ourType 2} 
&
{Weight 2 of Engel ourType 2} 
\\
\( 
\left[
\begin{array}{*{6}{c}}
\text{KerD} & 1 & 4 & 3 & 1 & 0 \\
\text{Bett} & 1 & 1 & 0 & 0 & 0 \\
\end{array}
\right]
\)
& 
\( 
\left[
\begin{array}{*{6}{c}}
\text{KerD} & 6 & 18 & 18 & 6 & 0 \\
\text{Bett} & 0 &  0 &  0 & 0 & 0 \\
\end{array}
\right]
\)
& 
\(\left[ 
\begin{array}{*{7}{c}}
\text{KerD} & 4 & 33 & 75 & 67 & 21 & 0 \\
\text{Bett} & 0 &  0 &  0 &  0 &  0 & 0 \\
\end{array}\right]
\)
\\
{Weight 0 of Engel ourType 3} 
& 
{Weight 1 of Engel ourType 3} 
&
{Weight 2 of Engel ourType 3} 
\\
\( 
\left[
\begin{array}{*{6}{c}}
\text{KerD} & 1 & 4 & 3 & 1 & 0 \\
\text{Bett} & 1 & 1 & 0 & 0 & 0 \\
\end{array}
\right]
\)
& 
\( 
\left[
\begin{array}{*{6}{c}}
\text{KerD} & 6 & 18 & 18 & 6 & 0 \\
\text{Bett} & 0 &  0 &  0 & 0 & 0 \\
\end{array}
\right]
\)
& 
\(\left[ 
\begin{array}{*{7}{c}}
\text{KerD} & 4 & 33 & 75 & 67 & 21 & 0 \\
\text{Bett} & 0 &  0 &  0 &  0 &  0 & 0 \\
\end{array}\right]
\)
\\
{Weight 0 of Engel ourType 4} 
& 
{Weight 1 of Engel ourType 4} 
&
{Weight 2 of Engel ourType 4} 
\\
\( 
\left[
\begin{array}{*{6}{c}}
\text{KerD} & 1 & 4 & 3 & 1 & 0 \\
\text{Bett} & 1 & 1 & 0 & 0 & 0 \\
\end{array}
\right]
\)
& 
\( 
\left[
\begin{array}{*{6}{c}}
\text{KerD} & 6 & 18 & 18 & 6 & 0 \\
\text{Bett} & 0 &  0 &  0 & 0 & 0 \\
\end{array}
\right]
\)
& 
\(\left[ 
\begin{array}{*{7}{c}}
\text{KerD} & 4 & 33 & 75 & 67 & 21 & 0 \\
\text{Bett} & 0 &  0 &  0 &  0 &  0 & 0 \\
\end{array}\right]
\)
\\
{Weight 0 of Engel ourType 5} 
& 
{Weight 1 of Engel ourType 5} 
&
{Weight 2 of Engel ourType 5} 
\\
\( 
\left[
\begin{array}{*{6}{c}}
\text{KerD} & 1 & 4 & 3 & 1 & 1 \\
\text{Bett} & 1 & 1 & 0 & 1 & 1 \\
\end{array}
\right]
\)
& 
\( 
\left[
\begin{array}{*{6}{c}}
\text{KerD} & 6 & 18 & 18 & 6 & 0 \\
\text{Bett} & 0 &  0 &  0 & 0 & 0 \\
\end{array}
\right]
\)
& 
\(\left[ 
\begin{array}{*{7}{c}}
\text{KerD} & 4 & 33 & 77 & 67 & 23 & 1 \\
\text{Bett} & 0 &  2 &  2 &  2 &  3 & 1 \\
\end{array}\right]
\)
\\
{Weight 0 of Engel ourType 6} 
& 
{Weight 1 of Engel ourType 6} 
&
{Weight 2 of Engel ourType 6} 
\\
\( 
\left[
\begin{array}{*{6}{c}}
\text{KerD} & 1 & 4 & 3 & 1 & 1 \\
\text{Bett} & 1 & 1 & 0 & 1 & 1 \\
\end{array}
\right]
\)
& 
\( 
\left[
\begin{array}{*{6}{c}}
\text{KerD} & 6 & 18 & 18 & 6 & 0 \\
\text{Bett} & 0 &  0 &  0 & 0 & 0 \\
\end{array}
\right]
\)
& 
\(\left[ 
\begin{array}{*{7}{c}}
\text{KerD} & 4 & 33 & 77 & 67 & 21 & 2 \\
\text{Bett} & 0 &  2 &  2 &  0 &  2 & 2 \\
\end{array}\right]
\)
\end{tabular}
%\end{center}
\caption{by Engel-mySols-v2.mpl}
\label{Engel-mySols}
\end{table}

\kmcomment{ 
\newpage 
% offset=-65mm -35mm, clip,  pagecommand={\thispagestyle{plain}} ]
\begin{minipage}{0.3\textwidth}
\includepdf[pages={-},scale=0.8, frame=false, trim = 70mm 0mm 70mm 0mm, 
offset=-65mm 25mm, clip,  pagecommand={\thispagestyle{plain}} ]
{Engel-mySols-wt0.pdf}
\end{minipage}
\begin{minipage}{0.35\textwidth}
\includepdf[pages={-},scale=0.8, frame=false, trim = 70mm 0mm 70mm 0mm, 
offset=-5mm 25mm, clip,  pagecommand={\thispagestyle{plain}} ]
{Engel-mySols-wt1.pdf}
\end{minipage}
\begin{minipage}{0.35\textwidth}
\includepdf[pages={-},scale=0.8, frame=false, trim = 70mm 0mm 70mm 0mm, 
offset=55mm 25mm, clip,  pagecommand={\thispagestyle{plain}} ]
{Engel-mySols-wt2.pdf}
\end{minipage}
}
%\newpage
\begin{comment}
\newpage 
% offset=-65mm -35mm, clip,  pagecommand={\thispagestyle{plain}} ]
\begin{minipage}{0.3\textwidth}
\includepdfmerge[nup=3x1, scale=0.85, frame=false, trim = 70mm 0mm 70mm 0mm, 
offset=0mm 20mm, clip,  pagecommand={\thispagestyle{plain}} ]
{Engel-mySols-wt0.pdf,Engel-mySols-wt1.pdf,Engel-mySols-wt2.pdf }
\end{minipage}

%\includepdf[pages={1}, frame=true]{Engel-1-012.pdf}

%$\hrulefill BEGIN \hrulefill $

%\includegraphics[scale=0.75, trim=0mm 105mm 0mm 10mm, clip]{Engel-1-012.pdf}
%\includepdf[pages={1-},scale=1.0, frame=false, 
% pagecommand={\thispagestyle{plain}} ]
% {July17-Engel-1-012.pdf}

\kmcomment{ 
\includepdf[pages={1},scale=1.0, trim=20mm 0mm 15mm 0mm, frame=false]{Engel-1-012.pdf}
\includepdf[pages={1},scale=1.0, trim=20mm 0mm 15mm 0mm, clip=true]{Engel-1-012.pdf}

\includepdf[pages={1},scale=1.0, trim=20mm -50mm 15mm 80mm, frame=false]{Engel-1-012.pdf}
\includepdf[pages={1},scale=1.0, trim=20mm -50mm 15mm 80mm, clip=true]{Engel-1-012.pdf}
}
%$\dotfill END \dotfill $

\kmcomment{ 
\includepdf[pages={1-},scale=1.0, frame=false, 
pagecommand={\thispagestyle{plain}} ]
{July17-Engel-2-012.pdf}

\includepdf[pages={1-},scale=1.0, frame=false,
pagecommand={\thispagestyle{plain}} ]
{July17-Engel-3-012.pdf}

\includepdf[pages={1-},scale=1.0, frame=false,
pagecommand={\thispagestyle{plain}} ]
{July17-Engel-4-012.pdf}
 }

\end{comment}

%\vspace*{125mm}
%\newpage

Then we have the proposition below.   

\begin{prop}
By the three weighted generic results about Betti numbers for each type \(
\texttt{Type(i)}\) where \( i=1\ldots 6 \),  those six types are  divided
into 4 classes Type(1), \{Type(2),  Type(3),  Type(4)\},   Type(4) and
Type(6), so that each two Lie algebras are not isomorphic if they belong
different classes.

\end{prop}

\begin{remark}
We expect to find more refinement if we try the third  or more weight cases.    
The word ``generic'' in the theorem above means we just follow Maple formal 
manipulations. When we want the kernel dimension, we have two ways: one is
solve the linear equations, and the number of free parameters means the
kernel dimension. The other is to fix the Groebner basis whose
number shows the rank.  
For instance, assume we have a linear system \(\{ a x + b y =0 \}\)
for \(x, y\) are independent variables and \(a,b\) are general constants. 
Maple !solve!-command shows  \( x = - \frac{b y}{a}, y=y\) and 
Maple !Groebner:-Basis({a x + b y}, tdeg( x,y ))! shows \([ a x + b y ]\)
formally. We say those answers are ``generic''.  
If \(a \ne 0\) or \( b\ne 0\), then the result is included in the generic solution.  
But, if \(a=b=0\), then the kernel space is 2-dimensional and the Groebner basis
is 0-dimensional. We may be careful Maple does not care about this case.  
\end{remark}

\kmcomment{ There are several questions. 
\newpage
\begin{minipage}{\textwidth}
\includepdf[pages=1 ]{cropOut.pdf}
\end{minipage}
%\includepdf[pages=2-]{cropOut.pdf}

%\includegraphics{eng-try-1-out.pdf}
%\includegraphics[]{cropOut.pdf}
%\includepdf[]{cropOut.pdf}
\includepdf{eng-try-1-out.pdf}

This space is Mottainai.

}

% \newpage
\section{Application of super homology of cotangent type
\(\ds \mathop{\oplus}_{\ell=0} ^{4} \Lambda ^{\ell} \frakg^{*}\) 
}

We have superalgebras on the exterior algebra of differential forms of
manifold .  We restrict our attention to the invariant forms of Lie groups,
and have the super homology groups.  We refer to \cite{Mik:Miz:superForms}
about precise definitions, but in this note, the grading for forms is sharp,
namely, the grade of \(p\)-form is \( -(p +1) \). On the contrary, the grade
of \(p\)-multi vector field is \(p-1\).  Given a (primary) weight \(wt\), if
\(wt\) is non-negative, then the chain space is empty. Otherwise we have
chain space bounded to \( |wt|\),  which is 1-dimensional, consists of \(\ds
\underbrace{1\mw \cdots \mw 1}_{|wt|}\), which is a kernel element.

\begin{wrapfigure}[5]{r}[5mm]{0.25\textwidth}
\vspace{-\baselineskip}
\(
\begin{array}{|c | *{5}{r}|}\hline
m & 1 & 2 & 3 & 4 & 5\\\hline
\text{SpaceDim} & 1 & 28 & 12 & 4 & 1 \\
\text{KerDim} & 1 &   &   &    &  1 \\
\text{Betti} &   &   &   &    &  1 \\\hline
\end{array}
\)
\end{wrapfigure}
\strut
We concentrate for the weight \(-5\) here. We follow the definition of \(m\)
-th chain space. It consists of  
\(\ds \frakhN{-1}^{a_{1}}\we \frakhN{-2}^{a_{2}}\we \frakhN{-3}^{a_{3}}\we \frakhN{-4}
^{a_{4}}\we \frakhN{-5}^{a_{5}}\) where \(\ds 
a_{1} + \cdots + a_{5} = m \) and \(\ds 
(-1) a_{1} + \cdots + (-5) a_{5} = -5 \). 
Directly, we see that \( m <= 5\) and \( a_{i} \) comes from the Young
diagrams of area 5 and length \(m\). 
Thus the common table becomes as
seen on the right end. 
When \(m=5\), the chain space is \( \frakhN{-1}^{5} \) and 1-dimensional. 
When \(m=4\), the chain space is \( \frakhN{-1}^{3} \mw \frakhN{-2} \) and 4-dimensional. 
When \(m=3\), the chain space is \(( \frakhN{-1} \mw \frakhN{-2}^{2} )\oplus
( \frakhN{-1}^{2} \mw \frakhN{-3}) 
\) and 12-dimensional. 
When \(m=2\), the chain space is \((\frakhN{-2} \mw \frakhN{-3})\oplus 
(\frakhN{-1} \mw \frakhN{-4} )
\) and 28-dimensional. 
When \(m=1\), the chain space is \( \frakhN{-5} = \verb|&^|( \zb{1},\zb{2},
\zb{3},\zb{4} )\) and 1-dimensional. 
By manipulating weight \((-5)\) super homology groups below, we classify 6
types into 3 as Type1, Type2 -- Type4, and Type5 -- Type 6, which is not
new in computation of the previous section. 

%\newpage

% \includepdf[pages={-}, pagecommand={}]{DualEngel-1-6-5.pdf}

\kmcomment{ Not Good 
%\begin{minipage}{\textwidth}
%\includepdf[pages=3, pagecommand={}]{DualEngel-1-6-5.pdf}
%\end{minipage}
\begin{minipage}{\textwidth}
\includepdf[pages=3, offset=0mm 80mm, trim=0mm 160mm 0mm 10mm, clip=true,
frame=true]{DualEngel-1-6-5.pdf}
\end{minipage}

%\vspace{-50mm}
ABC

}

%\newpage 

\begin{minipage}{\textwidth}
\includepdf[pages={-}, scale=0.9,offset = -10mm -155mm, trim = 15mm 0mm
45mm 0mm, frame=false,
pagecommand={\thispagestyle{plain}} ]{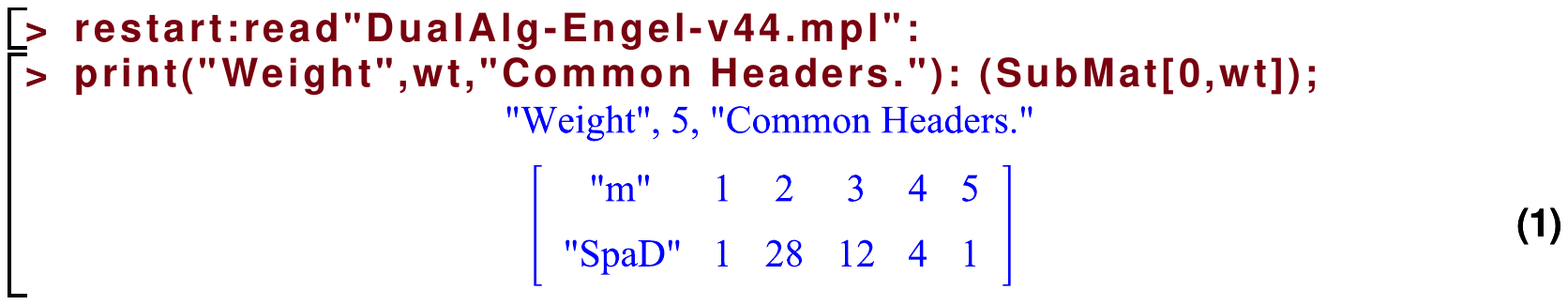}
\end{minipage}

%\newpage 

\begin{minipage}{\textwidth}
\includepdf[pages={-}, scale=0.9,offset = -10mm -190mm, trim = 15mm 0mm
45mm 0mm, frame=false,
pagecommand={\thispagestyle{plain}} ]{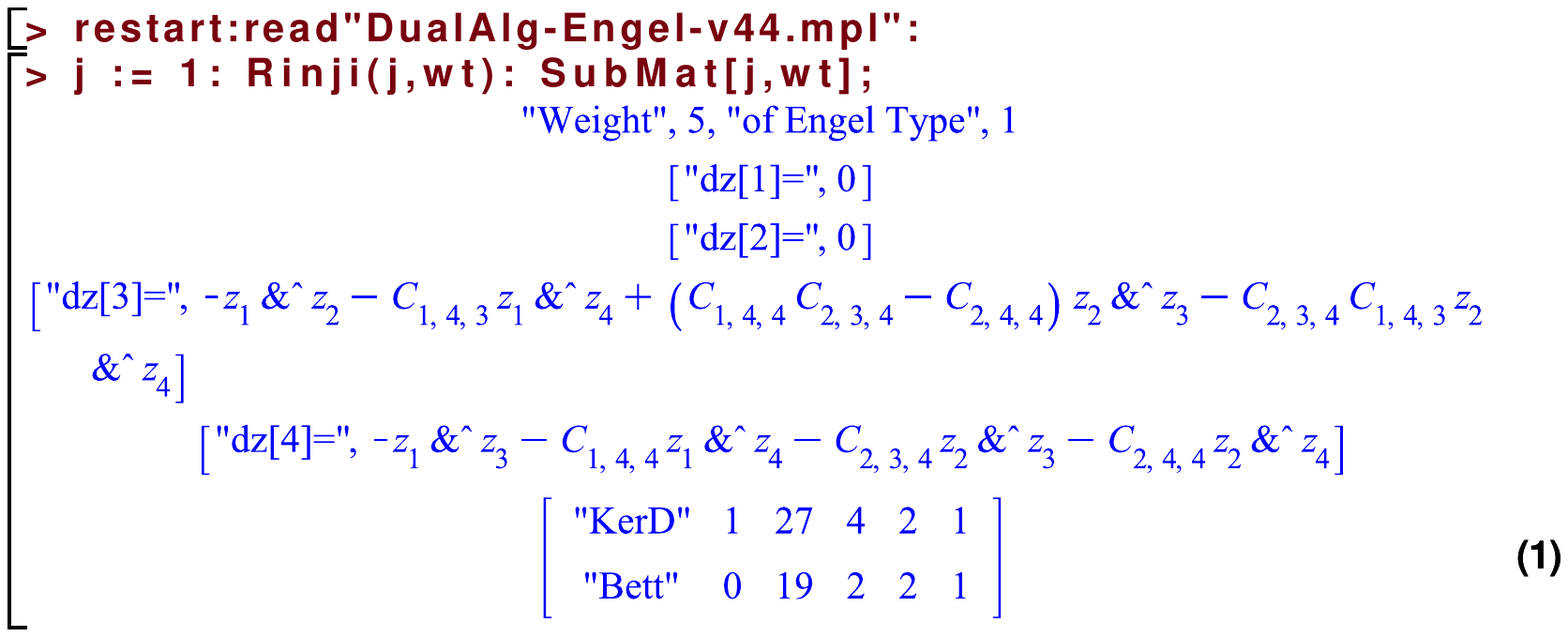}
\end{minipage}

\newpage

\begin{minipage}{\textwidth}
\includepdf[pages={-}, scale=0.9,offset = -10mm 10mm, trim = 15mm 0mm
45mm 0mm, frame=false,
pagecommand={\thispagestyle{plain}} ]{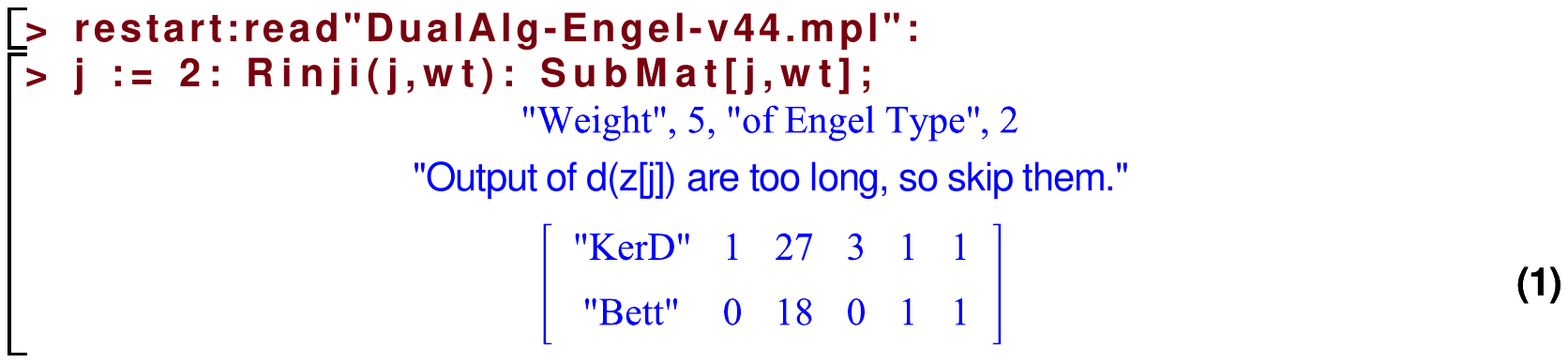}
\end{minipage}

\begin{minipage}{\textwidth}
\includepdf[pages={-}, scale=0.9,offset = -10mm -30mm, trim = 15mm 0mm
45mm 0mm, frame=false,
pagecommand={\thispagestyle{plain}} ]{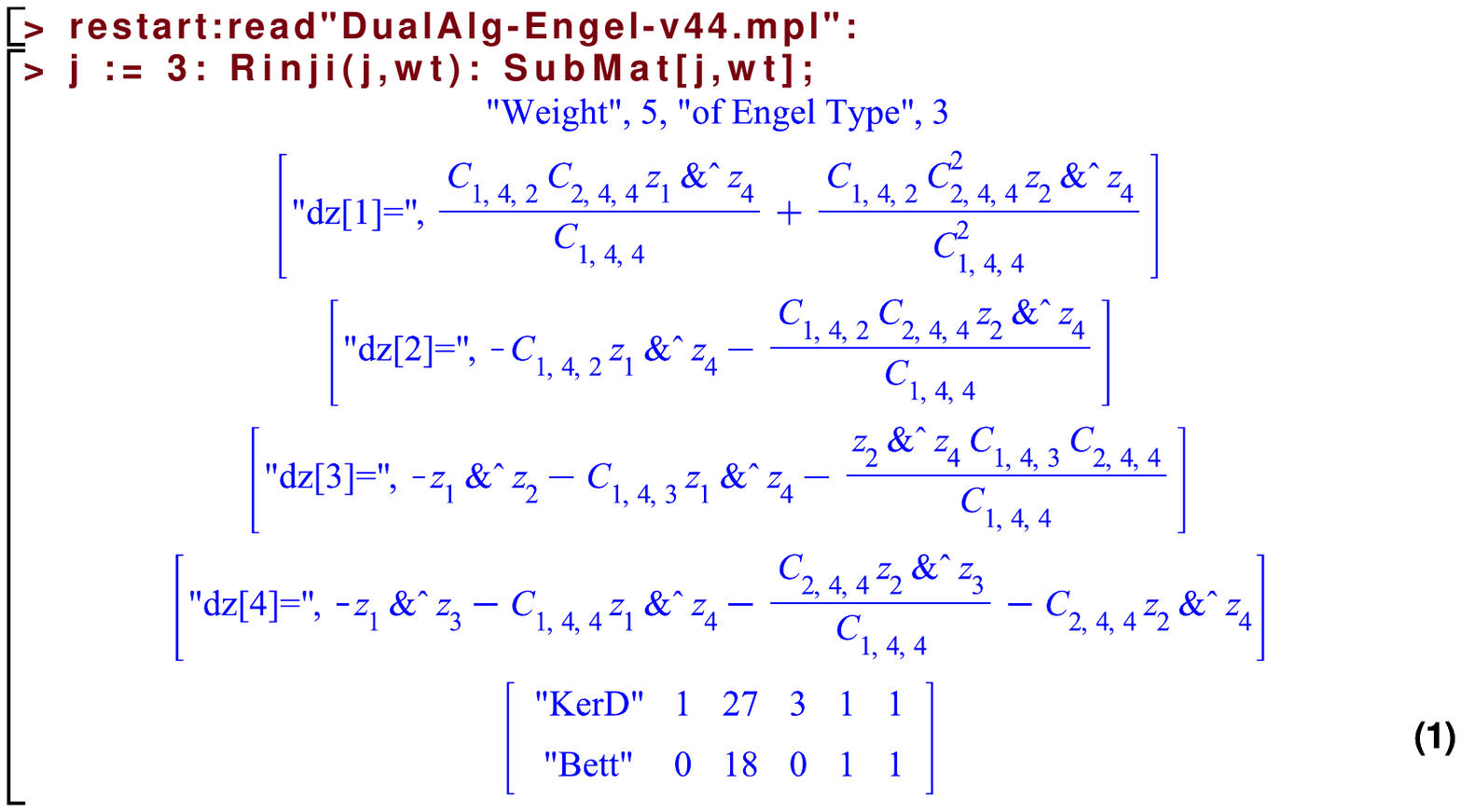}
\end{minipage}

\begin{minipage}{\textwidth}
\includepdf[pages={-}, scale=0.9,offset = -10mm -120mm, trim = 15mm 0mm
45mm 0mm, frame=false,
pagecommand={\thispagestyle{plain}} ]{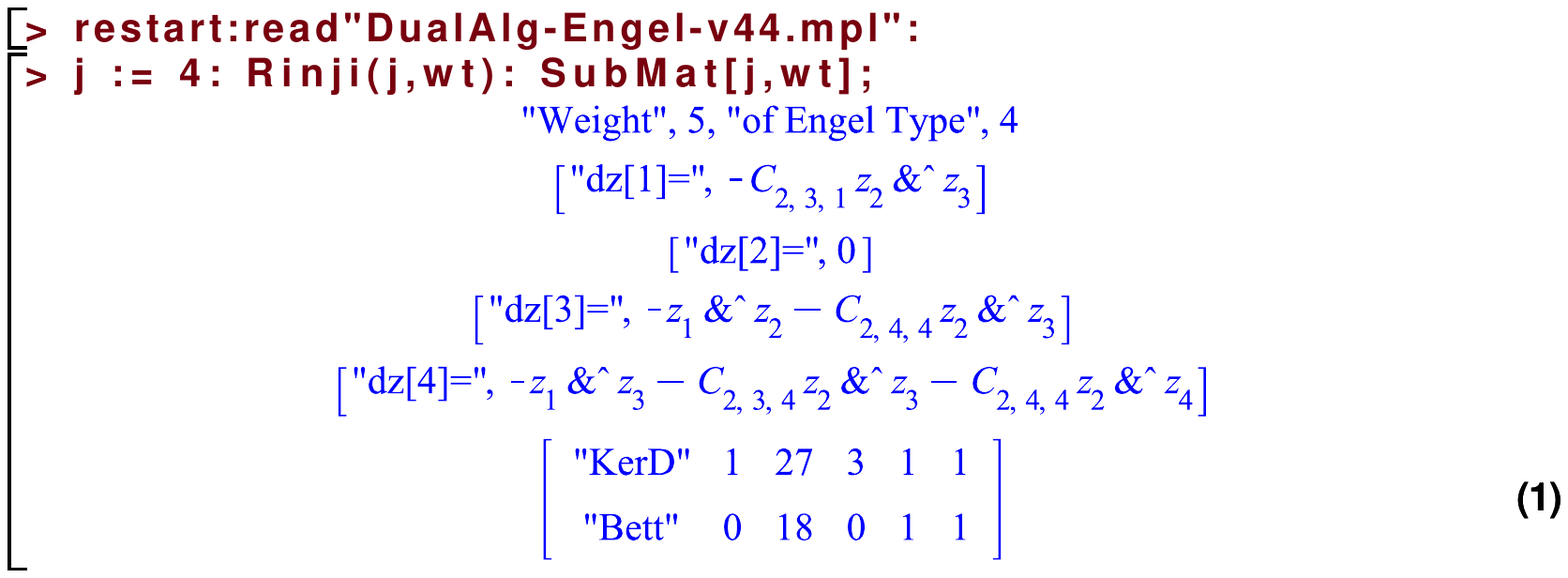}
\end{minipage}

\begin{minipage}{\textwidth}
\includepdf[pages={-}, scale=0.9,offset = -10mm -185mm, trim = 15mm 0mm
45mm 0mm, frame=false,
pagecommand={\thispagestyle{plain}} ]{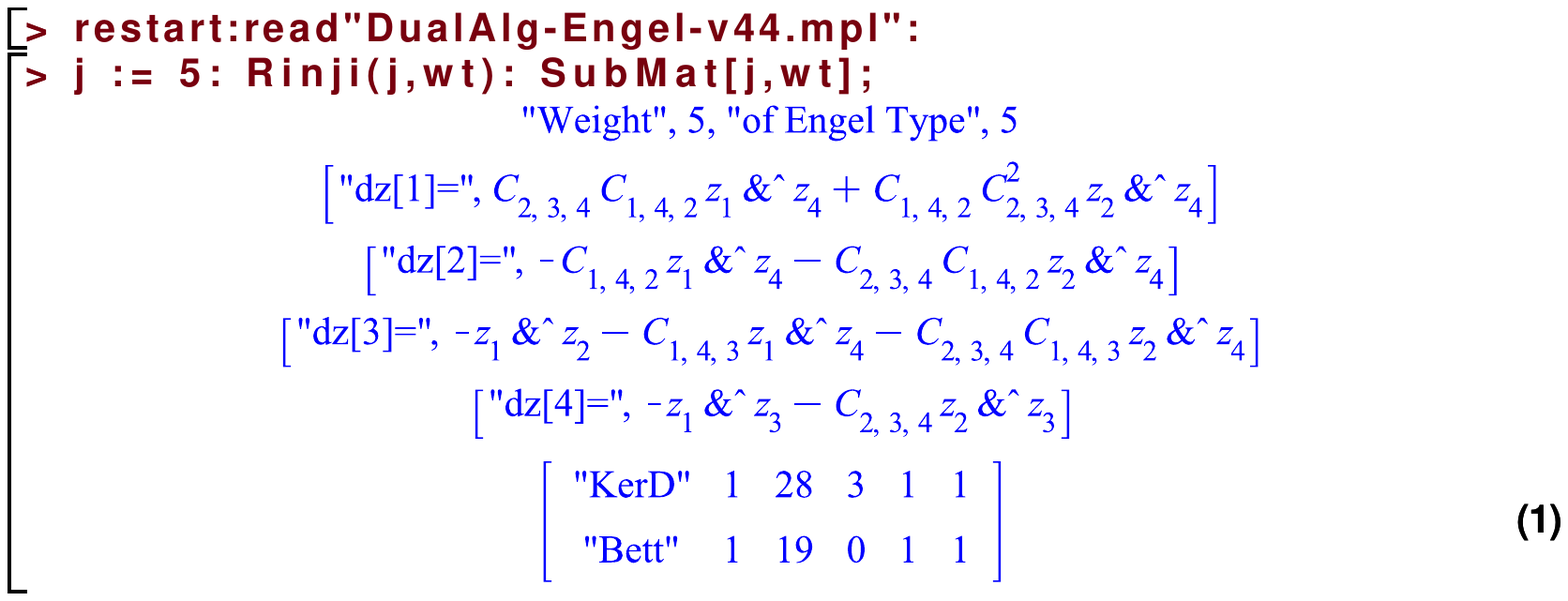}
\end{minipage}

\newpage

\begin{minipage}{\textwidth}
\includepdf[pages={-}, scale=0.9,offset = -10mm 0mm, trim = 15mm 0mm
45mm 0mm, frame=false,
pagecommand={\thispagestyle{plain}} ]{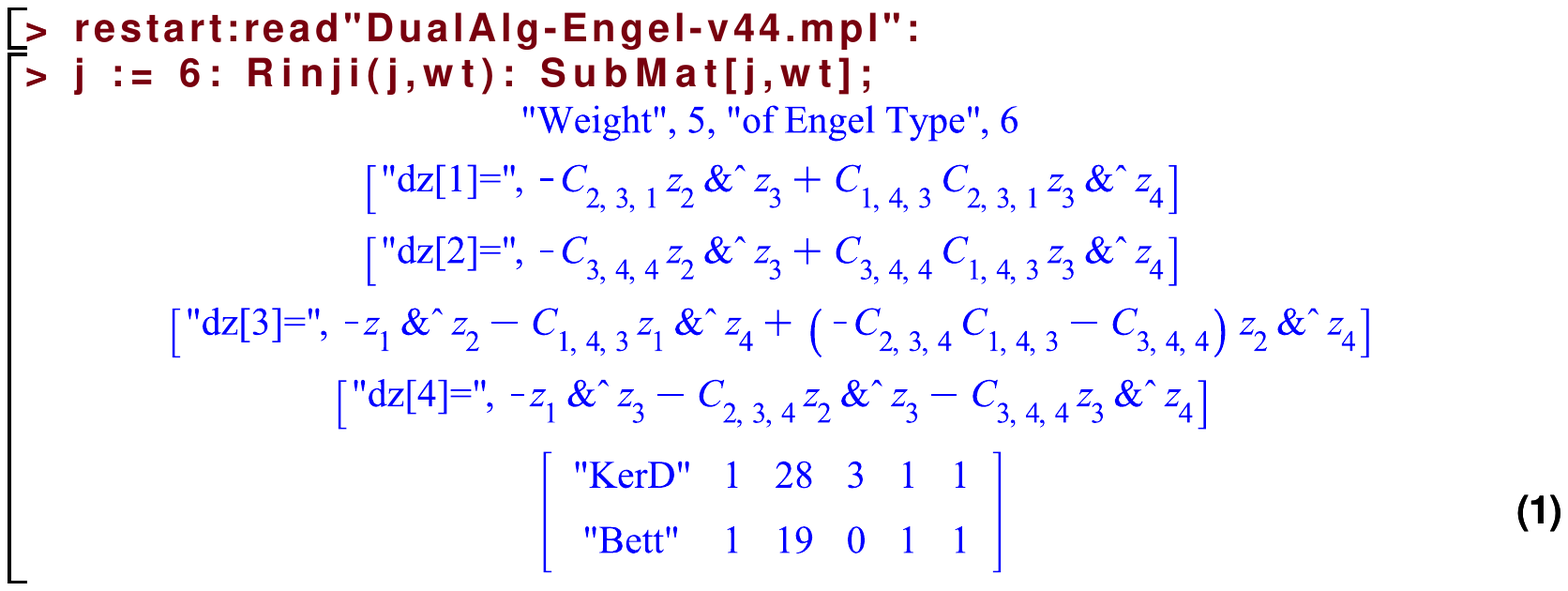}
\end{minipage}

\begin{comment}
\newpage 

\begin{minipage}{\textwidth}
\includepdf[nup=2x1, pages={-}, scale=0.9,offset = -10mm -10mm, trim = 15mm 0mm
45mm 0mm, frame=false,
pagecommand={\thispagestyle{plain}} ]{DualAlg-Engel-wt5-jul31.pdf}
\end{minipage}

\newpage 

\begin{minipage}{\textwidth}
\includepdf[pages={1}, scale=0.9,offset = 0mm -5mm, trim = 0mm 0mm 0mm
0mm, frame=false,
pagecommand={\thispagestyle{plain}} ]{DualEngel-1-6-5.pdf}
\end{minipage}

\includepdf[pages={2}, scale=0.9,offset = 0mm 5mm, trim = 0mm 0mm 0mm
0mm, frame=false,
pagecommand={\thispagestyle{plain}} ]{DualEngel-1-6-5.pdf}

\begin{minipage}{\textwidth}
\includepdf[pages={3}, trim = 0mm 100mm 0mm 0mm, offset=0mm 50mm, 
scale=0.9, clip , frame=false, 
pagecommand={\thispagestyle{plain}} ]{DualEngel-1-6-5.pdf}
\end{minipage}
\end{comment}

\vspace*{60mm}
\begin{remark}
Here, we explain about ``generic'' which was stated in the theorem 2.1 before, 
by outputs in this section more precisely.   
In general, we follow both simple !solve! and !Groebner:-Basis!. An easy
checkpoint is !denominator! or possibility of enumeration of elements of !Groebner:-Basis!.
\newcommand{\qb}[1]{q_{#1}}

% \texttt{Type(1)}: \vspace{-3mm}
\begin{entry}{7mm}
\item[Type(1): The case of \(m=2\):]    
\begin{align*}
SOL &= \qb{1} = \qb{1}, \ldots, \qb{11} = \qb{11}, \qb{13} = \qb{13} , \ldots,  \qb{28} = \qb{28};
\\ & 
\qb{12 } = - 
(\CC{2,3,4}\qb{6}-\CC{2,3,4}\qb{15}
+\CC{2,3,4}\qb{20}-\CC{2,3,4}\qb{27}-\qb{17}+\qb{22}
- \qb{28}\\&\qquad \quad  
- \frac{2}{\CC{1,4,4}}( -\CC{2,4,4}\qb{6}+\CC{2,4,4}\qb{15}-\CC{2,4,4}\qb{20}+\CC{2,4,4}\qb{27}) 
) 
\\
%\end{align*} \begin{align*}
EGB &= [\qb{6}(\CC{1,4,4}\CC{2,3,4}-2\CC{2,4,4})+\qb{12}\CC{1,4,4}+(-\CC{1,4,4}\CC{2,3,4}+
2\CC{2,4,4})\qb{15}-\qb{17}\CC{1,4,4}+\qb{20}(\CC{1,4,4}\CC{2,3,4}-2\CC{2,4,
4})\\
&  +\qb{22}\CC{1,4,4}+(-\CC{1,4,4}\CC{2,3,4}+2\CC{2,4,4})\qb{27}-\qb{28}\CC{1,4,4}]:
\end{align*}
\(\CC{1,4,4}\ne 0\) is in generic. 
If 
\(\CC{1,4,4}= 0\) then ORG = EGB = [\( -2 \CC{2,4,4} ( \qb{6}- \qb{15}+ \qb{20}
- \qb{27}) \) ], and so if \(\CC{2,4,4}\ne 0\) then the case is in generic. 
If 
\(\CC{1,4,4}= 0\) and \(  \CC{2,4,4}=0 \), then rank is 0 and the kernel dim is
28. 
\item[ Type(1): 
The case of \(m=3\):]
\begin{align*}
SOL &= \{\qb{1} = \qb{1}, \qb{2} = 0, \qb{3} = 0, \qb{4} = 0, \qb{5} = 0, \qb{6} = 0, \qb{7} = 
\qb {7}, \qb{8} = \qb{8}, \qb{9} = \qb{9},\\
& \qquad \qb{10} = (\CC{1,4,4}^2 \CC{2,3,4} \qb{8}+\CC{1,4,3} \CC{2,3
,4} \qb{8}-\CC{1,4,4} \CC{2,3,4} \qb{9}-\CC{1,4,4} \CC{2,4,4} \qb{8}+\CC{2,4,4} \qb{9})/\CC{1,4,3},
\\&\qquad 
\qb{11} = -\CC{1,4,4} \CC{2,3,4} \qb{8}+\CC{2,3,4} \qb{9}+\CC{2,4,4} \qb{8},
\qb{12} = 0\};
\\
EGB &= [\qb{12}, \qb{9} (\CC{1,4,3} \CC{2,3,4}^2+\CC{1,4,4} \CC{2,3,4} \CC{2,4,4}-\CC{2,4,4}^2)
+\qb{10} (-\CC{1,4,3} \CC{1,4,4} \CC{2,3,4}+\CC{1,4,3} \CC{2,4,4})
\\&\qquad\quad +\qb{11} (-\CC{1,4,4}^2 \CC{2,3 ,4}-\CC{1,4,3} \CC{2,3,4}+\CC{1,4,4} \CC{2,4,4}),\\&
(\CC{1,4,3} \CC{2,3,4}^2+\CC{1,4,4} \CC{2,3,4 } \CC{2,4,4}-\CC{2,4,4}^2) \qb{8}-\CC{1,4,3} \CC{2,3,4} \qb{10}
+(-\CC{1,4,4} \CC{2,3,4}+\CC{2,4,4 }) \qb{11}, \qb{6}, \qb{5}, \qb{4}, \qb{3}, \qb{2}];
\end{align*}
Thus, if \( \CC{1,4,3} \ne 0\) then it is generic and rank is 8.   
If \( \CC{1,4,3} =0 \),  then 
\(
EGB = [\qb{12}, A, B, \qb{6}, \qb{5}, \qb{4}, \qb{3}, \qb{2}]
\) where 
\begin{alignat*}{3} 
A &= 
\kmcomment{
\qb{9} 
( \CC{1,4,4} \CC{2,3,4} -\CC{2,4,4})\CC{2,4,4}
 +\qb{11}\CC{1,4,4} (-\CC{1,4,4} \CC{2,3 ,4}+ \CC{2,4,4})
 =   }
( \CC{1,4,4} \CC{2,3,4} -\CC{2,4,4})(  \CC{2,4,4} \qb{9} - \CC{1,4,4}\qb{11})
\;, 
&\quad\text{and}\quad & 
B &= \kmcomment{ 
(\CC{1,4,4} \CC{2,3,4 } \CC{2,4,4}-\CC{2,4,4}^2) \qb{8}
+(-\CC{1,4,4} \CC{2,3,4}+\CC{2,4,4 }) \qb{11}, 
=  }
(\CC{1,4,4} \CC{2,3,4 } -\CC{2,4,4})(\CC{2,4,4}\qb{8} - \qb{11})  \; .
\end{alignat*}
If \( \CC{1,4,4} \CC{2,3,4 } -\CC{2,4,4}=0 \), \(A=B=0\) and the rank is 6.  
If \( \CC{1,4,4} \CC{2,3,4 } -\CC{2,4,4}\ne 0 \), then  \(B\ne 0\) and  
\(A' = \CC{2,4,4} r_{9} -\CC{1,4,4} r_{11}\) 
\(B' = \CC{2,4,4} r_{8} - r_{11}\). The wedge product \( \&{}\hat( A',B' ) =
-\CC{2,4,4}( \CC{2,4,4} \&{}\hat(r_{8},r_{9})-  \CC{1,4,4} \&{}\hat(r_{8},r_{10}) +
\&{}\hat(r_{9},r_{10}) \).  
Thus, if \(\CC{2,4,4}=0\) then \(A = \CC{1,4,4} B\) and the rank is 7. 
If \(\CC{2,4,4}\ne 0\) then \(A\) and \(B\) are linearly independent and the
rank is 8 (in generic). 
\item[Type(1): Other cases are in generic.] 
\end{entry}
By the same discussion, we have checked the types 2, 3, 5 and 6 are only
generic. For the types 2 and 3, \(\CC{1,4,4}\ne 0\) is assumed.  

%\texttt{Type(4)}:
\vspace{-3mm}
\begin{entry}{7mm}
\item[Type(4): The case of \(m=2\):]    
The kernel condition is \([ -2 \CC{2,4,4}( \qb{6}- \qb{15} + \qb{20}-\qb{27})
]\). So if \(\CC{2,4,4} \ne 0\) then the kernel dimension is 27 and the rank
is 1 (in generic). 
If \(\CC{2,4,4} = 0\) then the kernel dimension is 28 and the rank is 0. 
\item[Type(4): Other cases are in generic.] 
\end{entry}
\end{remark}

\begin{table}
\centering
%\begin{center}
\begin{tabular}{c@{\hspace{20mm}}c}
{Weight 6, Common Headers} 
&
{Weight 7, Common Headers} 
\\
\( 
\left[
\begin{array}{*{6}{c}}
\text{m}    &  2 & 3 & 4 & 5 & 6 \\
\text{SpaD} & 38 & 32 & 12 & 4 & 1 \\
\end{array}
\right]
\)
& 
\(\left[ 
\begin{array}{*{7}{c}}
\text{m} &    2 & 3 & 4 & 5 & 6 & 7 \\
\text{SpaD} & 28 & 76 & 32 & 12 & 4 & 1 \\
\end{array}\right]
\)
\\
{Weight 6 of Engel ourType 1} 
&
{Weight 7 of Engel ourType 1} 
\\
\( 
\left[
\begin{array}{*{6}{c}}
\text{KerD} & 38 & 12 & 4 & 2 & 1 \\
\text{Bett} & 18 & 4 & 2 & 2 & 1 \\
\end{array}
\right]
\)
& 
\(\left[ 
\begin{array}{*{7}{c}}
\text{KerD} & 28 & 50 & 10 & 4 & 2 & 1 \\
\text{Bett} & 4 & 28 & 2 & 2 & 2 & 1 \\
\end{array}\right]
\)
\\
{Weight 6 of Engel ourType 2} 
&
{Weight 7 of Engel ourType 2} 
\\
\( 
\left[
\begin{array}{*{6}{c}}
\text{KerD} & 38 & 10 & 3 & 1 & 1 \\
\text{Bett} & 16 &  1 &  0 & 1 & 1 \\
\end{array}
\right]
\)
& 
\(\left[ 
\begin{array}{*{7}{c}}
\text{KerD} & 28 & 50 & 9 & 3 & 1 & 1 \\
\text{Bett} & 4 &  27 &  0 &  0 &  1 & 1 \\
\end{array}\right]
\)
\\
{Weight 6 of Engel ourType 3} 
&
{Weight 7 of Engel ourType 3} 
\\
\( 
\left[
\begin{array}{*{6}{c}}
\text{KerD} & 38 & 12 & 3 & 1 & 1 \\
\text{Bett} & 18 &  3 &  0 & 1 & 1 \\
\end{array}
\right]
\)
& 
\(\left[ 
\begin{array}{*{7}{c}}
\text{KerD} & 28 & 52 & 9 & 3 & 1 & 1 \\
\text{Bett} & 6 &  29 &  0 &  0 &  1 & 1 \\
\end{array}\right]
\)
\\
{Weight 6 of Engel ourType 4} 
&
{Weight 7 of Engel ourType 4} 
\\
\( 
\left[
\begin{array}{*{6}{c}}
\text{KerD} & 38 & 10 & 3 & 1 & 1 \\
\text{Bett} & 16 &  1 &  0 & 1 & 1 \\
\end{array}
\right]
\)
& 
\(\left[ 
\begin{array}{*{7}{c}}
\text{KerD} & 28 & 50 & 9 & 3 & 1 & 1 \\
\text{Bett} & 4 &  27 &  0 &  0 &  1 & 1 \\
\end{array}\right]
\)
\\
{Weight 6 of Engel ourType 5} 
&
{Weight 7 of Engel ourType 5} 
\\
\( 
\left[
\begin{array}{*{6}{c}}
\text{KerD} & 6 & 13 & 3 & 1 & 1 \\
\text{Bett} & 19 &  4 &  0 & 1 & 1 \\
\end{array}
\right]
\)
& 
\(\left[ 
\begin{array}{*{7}{c}}
\text{KerD} & 28 & 53 & 10 & 3 & 1 & 1 \\
\text{Bett} & 7 &  31 &  1 &  0 &  1 & 1 \\
\end{array}\right]
\)
\\
{Weight 6 of Engel ourType 6} 
&
{Weight 7 of Engel ourType 6} 
\\
\( 
\left[
\begin{array}{*{6}{c}}
\text{KerD} & 38 & 11 & 3 & 1 & 1 \\
\text{Bett} & 17 &  2 &  0 & 1 & 1 \\
\end{array}
\right]
\)
& 
\(\left[ 
\begin{array}{*{7}{c}}
\text{KerD} & 28 & 53 & 10 & 3 & 1 & 1 \\
\text{Bett} & 7 &  31 &  1 &  0 &  1 & 1 \\
\end{array}\right]
\)
\end{tabular}
%\end{center}
\caption{by DualAlg-Engel-v44.mpl}
\label{Dual-Engel}
\end{table}

%\end{document}

\begin{comment}

\newpage 
\begin{minipage}{0.5\textwidth}
%\includepdf [pages={-},scale=1.0, offset=-45mm -95mm, 
\includepdf [pages={-},scale=1.0, offset=-45mm -0mm, 
trim = 55mm 10mm 70mm 0mm,
clip, pagecommand={\pagestyle{plain}} ]{DualAlg-Engel-wt6-july29.pdf} 
\end{minipage}
%
\begin{minipage}{0.5\textwidth}
\includepdf [ pages={-},scale=1.0, offset=45mm -0mm, 
trim = 55mm 10mm 70mm 0mm,
clip, pagecommand={\pagestyle{plain}} ]{DualAlg-Engel-wt7-july29.pdf} 
\end{minipage}
%\includepdf[pages={2-},scale=0.6, offset=0mm 0mm , 
% pagecommand={\pagestyle{plain}} ] {DualAlg-Engel--7.pdf} 

\kmcomment{
\begin{tikzpicture}[remember picture, overlay]
    \node at (current page.center){ \includegraphics[page=1,scale=0.6]{DualAlg-Engel--6.pdf}}; 
\end{tikzpicture}
\includepdf[pages={2-last}] {DualAlg-Engel--6.pdf} 
  }

%\vspace*{100mm}
\end{comment}

%\newpage 
\section{Cases for the extended superalgebra 
\(\ds \mathop{\oplus}_{\ell=0} ^{4} \Lambda ^{\ell} \frakg^{*} \oplus \frakg \) 
% of 
% \(\ds \mathop{\oplus}_{\ell=0} ^{4} \Lambda ^{\ell} \frakg^{*}\) 
}
In 
\cite{Mik:Miz:superForms}, we know that   
\(\ds \frakg \oplus \mathop{\oplus}_{\ell=0} ^{4} \Lambda ^{\ell} \frakg^{*}\) 
is a super superalgebra of 
\(\ds \mathop{\oplus}_{\ell=0} ^{4} \Lambda ^{\ell} \frakg^{*}\) 
by using Lie derivation, in more general context. In this section, we try and
see the homology groups of those extended superalgebras for each 
Engel-like Lie algebra  \(\frakg\). The following output by Maple2021, which
are drove by ZeroPlus-Engel-v3.mpl and are shown below, implies 
that six types are divided into 5 classes. Right now, type 2 and type 4 have
the same table and are not distinguished on this job. However, using other
enhanced superalgebras with the weight \(-3\) in this section, we conclude those six types are not isomorphic as Lie
algebras. Thus, we claim the next theorem. 
\begin{thm}
The six types, which come from Lie algebras axioms, are mutually not
isomorphic in generic.  
\end{thm}
\begin{remark} When the weight is \(-2\), the output on the left hand side below says
nothing about Type 2 and Type 4. The weight \(-3\) output on the right hand
side claim that six types are mutually not isomorphic as Lie algebras. 
\end{remark}

\begin{table}
\centering
%\begin{center}
\begin{tabular}{c@{\hspace{20mm}}c}
{Weight 2, Common Headers} 
&
{Weight 3, Common Headers} 
\\
\( 
\left[
\begin{array}{*{7}{c}}
\text{m}   & 1 &  2 & 3 & 4 & 5 & 6 \\
\text{SpaD}& 4 & 17 & 28 & 22 & 8 & 1 \\
\end{array}
\right]
\)
& 
\(\left[ 
\begin{array}{*{8}{c}}
\text{m} & 1 &   2 & 3 & 4 & 5 & 6 & 7 \\
\text{SpaD} & 6 & 28 & 53 & 52 & 28 & 8 & 1 \\
\end{array}\right]
\)
\\
{Weight 2 of Engel ourType 1} 
&
{Weight 3 of Engel ourType 1} 
\\
\( 
\left[
\begin{array}{*{7}{c}}
\text{KerD} & 4 & 13 & 16 & 10 & 4 & 1 \\
\text{Bett} & 0 & 1 & 4 & 6 & 4 & 1 \\
\end{array}
\right]
\)
& 
\(\left[ 
\begin{array}{*{8}{c}}
\text{KerD} & 6 & 26 & 33 & 28 & 16 & 6 & 1 \\
\text{Bett} & 4 & 6 & 9 & 16 & 14 & 6 & 1 \\
\end{array}\right]
\)
\\
{Weight 2 of Engel ourType 2} 
&
{Weight 3 of Engel ourType 2} 
\\
\( 
\left[
\begin{array}{*{7}{c}}
\text{KerD} & 4 & 13 & 17 & 10 & 4 & 1 \\
\text{Bett} & 0 & 2 &  5 &  6 & 4 & 1 \\
\end{array}
\right]
\)
& 
\(\left[ 
\begin{array}{*{8}{c}}
\text{KerD} & 6 & 25 & 28 & 25 & 13 & 5 & 1 \\
\text{Bett} & 3 & 0 &  1 &  10 & 10 & 5 & 1 \\
\end{array}\right]
\)
\\
{Weight 2 of Engel ourType 3} 
&
{Weight 3 of Engel ourType 3} 
\\
\( 
\left[
\begin{array}{*{7}{c}}
\text{KerD} & 4 & 13 & 18 & 10 & 4 & 1 \\
\text{Bett} & 0 &  3 &  6 &  6 & 4 & 1 \\
\end{array}
\right]
\)
& 
\(\left[ 
\begin{array}{*{8}{c}}
\text{KerD} & 6 & 25 & 30 & 25 & 13 & 5 & 1 \\
\text{Bett} & 3 & 2 &  3 &  10 & 10 & 5 & 1 \\
\end{array}\right]
\)
\\
{Weight 2 of Engel ourType 4} 
&
{Weight 3 of Engel ourType 4} 
\\
\( 
\left[
\begin{array}{*{7}{c}}
\text{KerD} & 4 & 13 & 17 & 10 & 4 & 1 \\
\text{Bett} & 0 &  2 &  5 &  6 & 4 & 1 \\
\end{array}
\right]
\)
& 
\(\left[ 
\begin{array}{*{8}{c}}
\text{KerD} & 6 & 25 & 29 & 25 & 13 & 5 & 1 \\
\text{Bett} & 3 & 1 &  2 &  10 & 10 & 5 & 1 \\
\end{array}\right]
\)
\\
{Weight 2 of Engel ourType 5} 
&
{Weight 3 of Engel ourType 5} 
\\
\( 
\left[
\begin{array}{*{7}{c}}
\text{KerD} & 4 & 13 & 18 & 10 & 5 & 1 \\
\text{Bett} & 0 &  3 &  6 &  7 & 5 & 1 \\
\end{array}
\right]
\)
& 
\(\left[ 
\begin{array}{*{8}{c}}
\text{KerD} & 6 & 25 & 30 & 25 & 16 & 5 & 1 \\
\text{Bett} & 3 &  2 &  3 & 13 & 13 & 5 & 1 \\
\end{array}\right]
\)
\\
{Weight 2 of Engel ourType 6} 
&
{Weight 3 of Engel ourType 6} 
\\
\( 
\left[
\begin{array}{*{7}{c}}
\text{KerD} & 4 & 14 & 16 & 10 & 5 & 1 \\
\text{Bett} & 1 & 2 &  4 &  7 & 5 & 1 \\
\end{array}
\right]
\)
& 
\(\left[ 
\begin{array}{*{8}{c}}
\text{KerD} & 6 & 25 & 29 & 25 & 16 & 5 & 1 \\
\text{Bett} & 3 &  1 &  2 & 13 & 13 & 5 & 1 \\
\end{array}\right]
\)
\end{tabular}
%\end{center}
\caption{by ZeroPlus-Engel-v3.mpl}
\label{Fig:ZeroPlus}
\end{table}

%\end{document}

%\newpage
\kmcomment{ 
\begin{minipage}{0.5\textwidth}
\includepdf [ pages={-},scale=1, offset=-35mm -80mm, trim = 60mm 10mm 50mm 0mm,
clip, pagecommand={\pagestyle{plain}}]{ZeroPlus-Engel-wt2-july29.pdf} 
\end{minipage}
\begin{minipage}{0.5\textwidth}
\includepdf [ pages={-},scale=1, offset=55mm -80mm, trim = 60mm 10mm 50mm 0mm,
clip, pagecommand={\pagestyle{plain}}]{ZeroPlus-Engel-wt3-july29.pdf} 
\end{minipage}

\begin{minipage}{\textwidth}
\includepdfmerge [nup=2x1, delta= -20mm 0mm, scale=1, offset=15mm -80mm, trim = 60mm 10mm 50mm 0mm,
clip, pagecommand={\pagestyle{plain}}]{ZeroPlus-Engel-wt2-july29.pdf, ZeroPlus-Engel-wt3-july29.pdf} 
\end{minipage}

\newpage
}

\section{Characteristic foliations}
There is a notion of 
characteristic foliations in Engel theory, which is rank 1 distribution
\(\mathfrak{L}\) of Engel distribution \( D \) satisfying 
\( \Sbt{\mathfrak{L}}{D^{2}} \subset D^{2}\) where 
\(D^{2} := D + \Sbt{D}{D} \) (, and 
\(D^{3} := D^{2} + \Sbt{D^{2}}{D^{2}} \)). 
In our six cases, they are given by
\( a \yb{1} + b \yb{2}\) for \texttt{mySol[3]},  and 
\( \CC{2,3,4} \yb{1} - \yb{2}\) for the others,  up to constant. 
It will be interesting to study their contributions to super homology
discussion. 

\section{Another approach of finding Engel-like structures}
So far, we studied Engel-like structures first fix dimensional conditions 
and next check Lie algebra structure.  
As \cite{MR3030954} points out, there is a complete classification list of
4-dimensional Lie algebra by  \cite{MR0155871}.  
Since we 
unfortunately do not have access to the original paper \cite{MR0155871}, 
we use the list in \cite{MR3030954}. 
\newcommand{\uType}{\text{Type}}
\begin{align*}
\uType[1] &=\{ 
            [\yb{2},\yb{4}] = \yb{1}, 
            [\yb{3},\yb{4}] = \yb{2} \}; 
\quad 
\uType[2] =\{
            [\yb{1},\yb{4}] =  a   \yb{1},  
            [\yb{2},\yb{4}] = \yb{2}, 
            [\yb{3},\yb{4}] = \yb{2} + \yb{3} \}; 
\\
\uType[3]&=\{
            [\yb{1},\yb{4}] = \yb{1}, 
            [\yb{3},\yb{4}] = \yb{2} \}; 
\quad   
\uType[4] =\{
            [\yb{1},\yb{4}] = \yb{1}, 
            [\yb{2},\yb{4}] = \yb{1}+\yb{2}, 
            [\yb{3},\yb{4}] = \yb{2}+\yb{3} \}; 
            \\
\uType[5]&=\{ 
            [\yb{1},\yb{4}] = \yb{1}, 
            [\yb{2},\yb{4}] = a \yb{2}, 
            [\yb{3},\yb{4}] = b \yb{3}, \text{ where }  a b \ne  0 \};  
\\
\uType[6]&=\{ 
            [\yb{1},\yb{4}] = a \yb{1}, 
            [\yb{2},\yb{4}] = b \yb{2} - \yb{3}, 
[\yb{3},\yb{4}] = \yb{2} + b \yb{3}, \text{ where } a \ne 0 \text{ and }  b \geq 0\};  
\\
\uType[7]&=\{ 
            [\yb{1},\yb{4}] = 2 \yb{1}, 
            [\yb{2},\yb{3}] = \yb{1}, 
            [\yb{2},\yb{4}] = \yb{2} , 
            [\yb{3},\yb{4}] = \yb{2} + \yb{3} \}; 
\\
\uType[8]&=\{ 
            [\yb{2},\yb{3}] = \yb{1}, 
            [\yb{2},\yb{4}] = \yb{2} , 
            [\yb{3},\yb{4}] = - \yb{3} \}; 
\\
\uType[9]&=\{ 
            [\yb{1},\yb{4}] = (1+b)  \yb{1}, 
            [\yb{2},\yb{3}] = \yb{1}, 
            [\yb{2},\yb{4}] = \yb{2} , 
            [\yb{3},\yb{4}] = b  \yb{3}, \text{ where }   -1 < b \leqq 1   \};  
\\
\uType[10]&=\{
            [\yb{2},\yb{3}] = \yb{1}, 
            [\yb{2},\yb{4}] = -\yb{3} , 
            [\yb{3},\yb{4}] =  \yb{2} \}; 
\\
\uType[11]&=\{
            [\yb{1},\yb{4}] = 2 a \yb{1}, 
            [\yb{2},\yb{3}] = \yb{1}, 
            [\yb{2},\yb{4}] = a \yb{2}-\yb{3} , 
            [\yb{3},\yb{4}] =  \yb{2} + a \yb{3} \}; 
\\
\uType[12]&=\{
            [\yb{1},\yb{3}] = \yb{1}, 
            [\yb{1},\yb{4}] = -\yb{2}, 
            [\yb{2},\yb{3}] = \yb{2}, 
            [\yb{2},\yb{4}] = \yb{1}  
            \};
\end{align*}

In this section, we first fix each 4-dimensional Lie algebra and try to find
2-dimensional subspace \(D\) satisfying Engel-like structure.  We prepare a
small proposition. 

\begin{prop}
Let \(\frakg\) be a 4-dimensional Lie algebra with bracket relations by a
basis \((\yb{1}, \ldots, \yb{4})\).  

Take a 2-dimensional plane \(D\) where 
\( \wb{1} = \sum_{i=1} ^{4}  \ANYb{p}{i} \yb{i}\) 
and 
\( \wb{2} = \sum_{i=1} ^{4}  \ANYb{q}{i} \yb{i}\) are its basis. 
Put 
\( \wb{3} = \Akt{ \wb{1} }{\wb{2}} \)
and 
\( \wb{4} = \Akt{ \wb{1} }{ \wb{3}} \).  
%\( \wb{4} = \Akt{ \wb{1} }{\Akt{ \wb{1} }{\wb{2}} } \).  

If \( \{ \wb{1}, \wb{2}, \wb{3} \}\) and 
 \( \{ \wb{1}, \wb{2}, \wb{3}, \wb{4} \}\) are linearly independent, 
i.e.,  
\(  \wb{1} \we  \wb{2} \we  \wb{3} \ne 0\) and 
 \( \wb{1} \we \wb{2} \we \wb{3} \we \wb{4} \ne 0\)  
 for some \( \{\ANYb{p}{i}\}_{i=1}^{4},  
 \{\ANYb{q}{i}\}_{i=1}^{4}  \), 
 then \( (\frakg, D) \)  is an Engel-like structure.   
\end{prop}

In this note, we call the scalar \( (\wb{1} \we \wb{2} \we \wb{3} \we \wb{4})
\slash (\yb{1} \we \yb{2} \we \yb{3} \we \yb{4})\) above by the
\emph{E-l-C} (Engel-like coefficient).   
\newcommand{\Det}{\text{Det}}
Using the notation \( \Det(i,j) = 
\ANYb{p}{i}  \ANYb{q}{j} 
- 
\ANYb{p}{j}  \ANYb{q}{i} 
= 
\begin{vmatrix} 
\ANYb{p}{i} &  \ANYb{p}{j} \\ 
\ANYb{q}{i} &  \ANYb{q}{j} \\ 
\end{vmatrix}\), 
the Engel-like coefficient of Lie algebras \uType[i] (i=1..12) are given as
follows. 

\begin{align*}
\text{E-l-C of } \uType[1] &= \ANYb{p}{4} \Det(3,4)^{3}
\qquad \qquad \qquad \qquad 
\text{E-l-C of } \uType[2] = (a-1)^{2} \ANYb{p}{4} 
\Det(1,4) \Det(3,4)^{2}
\\
\text{E-l-C of } \uType[3] &=  \ANYb{p}{4} 
\Det(1,4) \Det(3,4)^{2}
\hspace{17mm}
\text{E-l-C of } \uType[4] =  \ANYb{p}{4} \Det(3,4)^{3}
\\
\text{E-l-C of } \uType[5] &= (a-1)(b-1)(a-b) \ANYb{p}{4} \Det(1,4) \Det(2,4) \Det(3,4)
\\
\text{E-l-C of } \uType[6] &= ((a-b)^{2}+1) \ANYb{p}{4} \Det(1,4) 
(  \Det(2,4)^2 +  \Det(3,4)^ 2) 
\\
\text{E-l-C of } \uType[7] &=   \Det(3,4)^{2} 
( \ANYb{p}{4} \Det(1,4) + \ANYb{p}{4} \Det(2,3) 
+ \ANYb{p}{3} \Det(3,4) )  
\\
\text{E-l-C of } \uType[8] &=  - 2  \Det(2,4)\Det(3,4)
( \ANYb{p}{4} \Det(1,4) - \ANYb{p}{3} \Det(2,4) 
- \ANYb{p}{2} \Det(3,4) )  
\\
\text{E-l-C of } \uType[9] &=  - (b-1) \Det(2,4)\Det(3,4)
\left( \ANYb{p}{3} \Det(1,4) 
+ b ( \ANYb{p}{4} \Det(1,4) - \ANYb{p}{2} \Det(3,4)) \right)  
\\
\text{E-l-C of} \uType[10] &=  ( \Det(2,4)^2 + \Det(3,4)^2 ) 
( \ANYb{p}{4} \Det(1,4) 
+  \ANYb{p}{2} \Det(2,4) + \ANYb{p}{3} \Det(3,4)) 
\\
\text{E-l-C of } \uType[11] &=  ( \Det(2,4)^2 + \Det(3,4)^2 ) 
( a^{2} \ANYb{p}{4} \Det(1,4) 
+ a \ANYb{p}{4} \Det(2,3) + \ANYb{p}{4} \Det(1,4)  + \ANYb{p}{2} \Det(2,4) 
 + \ANYb{p}{3} \Det(3,4) ) 
\\
\text{E-l-C of } \uType[12] &= \ANYb{p}{4} \Det(3,4) 
\left( \Det(1,3)^2 + \Det(1,4)^2 + \Det(2,3)^2 + \Det(2,4)^2 
+ 2 \Det(1,2)  \Det(3,4) 
\right) 
\end{align*}

In the case of \uType[1], 
taking \( p = [0,0,0,1]\) and \( q = [0,0,1,0]\), i.e.,  
\( \wb{1} = \yb{4}\) and \( \wb{2} = \yb{3}\) give an Engel-like structure. 

In the case of \uType[2], 
if \(a=1\), then there is no Engel-like structure. 
Assume \( a \ne 1\), then 
\( p = [0,0,0,1]\) and \( q = [1,0,1,0]\)   
give an Engel-like structure. 

In the case of \uType[3],  
\( p = [0,0,0,1]\) and \( q = [1,0,1,0]\)   
give an Engel-like structure like as the second half of \uType[2]. 

The case \uType[4] is the same with \uType[1].   

In the case of \uType[5], 
if \((a-1)(b-1)(a-b) =0\), then there is no Engel-like structure. 
Otherwise, 
\( p = [0,0,0,1]\) and \( q = [1,0,1,0]\)   
give an Engel-like structure. 

In the case of \uType[6],  
\( p = [0,0,0,1]\) and \( q = [1,1,1,0]\)   
give an Engel-like structure. 

In the case of \uType[7],  
\( p = [0,0,1,1]\) and \( q = [0,0,0,1]\)   
give an Engel-like structure. 

In the case of \uType[8],  
\( p = [0,0,0,1]\) and \( q = [1,1,1,0]\)   
give an Engel-like structure. 

In the case of \uType[9],  
if \( b-1 = 0\) then there is no Engel-like structure. 
Otherwise, 
\( p = [1,1,1,1]\) and \( q = [0,0,0,1]\)   
give an Engel-like structure. 

In the case of \uType[10],  
\( p = [0,0,0,1]\) and \( q = [1,0,1,0]\)   
give an Engel-like structure. 

In the case of \uType[11],  
\( p = [0,0,1,0]\) and \( q = [0,0,0,1]\)   
give an Engel-like structure. 

In the case of \uType[12],  
\( p = [0,1,0,1]\) and \( q = [0,1,1,0]\)   
give an Engel-like structure.

%\nocite{Igles:Marrero} 
%\nocite{MR2284790} 
\nocite{Mik:Miz:super1} 
\nocite{Mik:Miz:super2} 
\nocite{Mik:Miz:super3} 
\nocite{Mik:Miz:superDefoForms} 
%\nocite{MR3030954} % Ghanam - Thompson  
%\nocite{MR0155871} % Mubarakzyanov 
\nocite{MR3030673} % Adachi 
\nocite{MR2336679} % Adachi 

\nocite{MR3184730}
%AUTHOR = {\v{S}nobl, Libor and Winternitz, Pavel},

\nocite{MR3487553}
%AUTHOR = {Biggs, Rory and Remsing, Claudiu C.},

% \begin{spacing}{1.5}
\bibliographystyle{plain}
\bibliography{km_refs}

\appendix
%\begin{spacing}{0.80}
\section{Table of bracket for Engel-like Lie algebras }
\renewcommand{\Bkt}{\mathop{Bkt}}
\newcommand{\newBkt}[2]{[#1,#2]}
%%  mySol[1  ]%%%%%%%%%%%%%%%%%%%%%%%%%%%%%%%%%%%%%%%%%%%%%%
\textbf{ourSolDep[1]}
\begin{alignat*}{2}
% \newBkt{\yb{1}}{\yb{2}}&=     \yb{3}\\
% \newBkt{\yb{1}}{\yb{3}}&=     \yb{4}\\
\newBkt{\yb{1}}{\yb{4}}&=\CC{1,4,3}    \yb{3}+\CC{1,4,4}    \yb{4} & \qquad 
\newBkt{\yb{2}}{\yb{3}}&=(-\CC{1,4,4} \CC{2,3,4}+\CC{2,4,4})    \yb{3}+\CC{2,3,4}   \yb{4}\\
\newBkt{\yb{2}}{\yb{4}}&=(\CC{1,4,3} \CC{2,3,4})    \yb{3}+\CC{2,4,4}
\yb{4}& \qquad 
\newBkt{\yb{3}}{\yb{4}}&=0
\end{alignat*}
%% mySol[2] %%%%%%%%%%%%%%%%%%%%%%%%%%%%%%%%%%%%%%%%%%%%%%

\kmcomment{ } 
\textbf{ourSolDep[2]}
\begin{align*}
% \newBkt{\yb{1}}{\yb{2}}=    \yb{3},
% \newBkt{\yb{1}}{\yb{3}}=     \yb{4},
\newBkt{\yb{1}}{\yb{4}}&=\frac{\CC{1,4,4}^2+4 \CC{1,4,3}}{8}\left( 
- (\CC{1,4,4} \CC{2,3,4}-2 \CC{2,4,4}) \yb{1}    
-\CC{1,4,4} \yb{2}\right) +\CC{1,4,3} \yb{3}+\CC{1,4,4} \yb{4}
\\ 
\newBkt{\yb{2}}{\yb{3}} &=
\frac{
(\CC{1,4,4}^2 + 4 \CC{1,4,3})( \CC{1,4,4}\CC{2,3,4}-\CC{2,4,4} )
}
{ 2 \CC{1,4,4}^{2}} ( -  (\CC{1,4,4}\CC{2,3,4} - 2 \CC{2,4,4}) \yb{1} -
\CC{1,4,4} \yb{2} )
\\& 
+(-\CC{1,4,4} \CC{2,3,4}+\CC{2,4,4}) \yb{3}+\CC{2,3,4} \yb{4} \\
\newBkt{\yb{2}}{\yb{4}}&=-\frac{1}{8} (\CC{1,4,4} \CC{2,3,4}-2 \CC{2,4,4})
(\CC{1,4, 4}^2+4 \CC{1,4,3})\CC{2,3,4} \yb{1}
-\frac{1}{8} (\CC{1,4,4}^2+4 \CC{1,4,3}) \CC{1, 4,4} \CC{2,3,4} \yb{2}
\\ & - \frac{1}{2 \CC{1,4,4}}(\CC{1,4,4}^{3}\CC{2,3,4} 
+ 2 \CC{1,4,3}\CC{1,4,4} \CC{2,3,4}- \CC{1,4,4}^{2} \CC{2,4,4}- 4 \CC{1,4,3}
\CC{2,4,4})\yb{3}+\CC{2,4,4}\yb{4} \\ 
\newBkt{\yb{3}}{\yb{4}}&=
\frac{(\CC{1,4,4}^2 + 4 \CC{1,4,3})( \CC{1,4,4}\CC{2,3,4} - \CC{2,4,4} ) }
{ 8 \CC{1,4,4}^{2} }
\left ( 
(\CC{1,4,4}\CC{2,3,4}-\CC{2,4,4}) (\CC{1,4,4}^2 + 4 \CC{1,4,3}) \yb{1}
\right.
\\ & \left.  
+ \CC{1,4,4} (\CC{1,4,4}^2 + 4 \CC{1,4,3}) \yb{2}
+ 2 \CC{1,4,4}^{2} \yb{3}
+ 4 \CC{1,4,4} \yb{4}\right)
\end{align*}
%% mySol[3] %%%%%%%%%%%%%%%%%%%%%%%%%%%%%%%%%%%%%%%%%%%%%%
\begin{alignat*}{2}
\shortintertext{ \textbf{ourSolDep[3]} }
% \newBkt{\yb{1}}{\yb{2}}=     \yb{3},
% \newBkt{\yb{1}}{\yb{3}}=     \yb{4},
\newBkt{\yb{1}}{\yb{4}}&=
- \frac{  \CC{1,4,2} \CC{2,4,4} } { \CC{1,4,4} }
\yb{1}+\CC{1,4,2} \yb{2}+\CC{1,4,3}\yb{3}
+\CC{1,4,4}\yb{4} & \qquad 
\newBkt{\yb{2}}{\yb{3}}& = \frac{\CC{2,4,4}}{\CC{1,4,4}}\yb{4} \\
\newBkt{\yb{2}}{\yb{4}}& =
\frac{ \CC{2,4,4} }{ \CC{1,4,4}^{2} } \left(  
- \CC{1,4,2}\CC{2,4,4}  \yb{1} + \CC{1,4,2} \CC{1,4,4} \yb{2}
+ \CC{1,4,3} \CC{1,4,4} \yb{3} + \CC{1,4,4}^{2}\yb{4}\right) & \qquad 
\newBkt{\yb{3}}{\yb{4}}&=0
\end{alignat*}
%% mySol[4] %%%%%%%%%%%%%%%%%%%%%%%%%%%%%%%%%%%%%%%%%%%%%%
\begin{alignat*}{2}
\shortintertext{ \textbf{ourSolDep[4]} }
% \newBkt{\yb{1}}{\yb{2}}=     \yb{3},
% \newBkt{\yb{1}}{\yb{3}}=     \yb{4},
\newBkt{\yb{1}}{\yb{4}}& =0 & \qquad 
\newBkt{\yb{2}}{\yb{3}}& =\CC{2,3,1} \yb{1}+\CC{2,4,4}\yb{3}+\CC{2,3,4} \yb{4} \\
\newBkt{\yb{2}}{\yb{4}}& =\CC{2,4,4}    \yb{4} & \qquad 
\newBkt{\yb{3}}{\yb{4}}&=0 
\end{alignat*}
%% mySol[5] %%%%%%%%%%%%%%%%%%%%%%%%%%%%%%%%%%%%%%%%%%%%%%
\begin{alignat*}{2}
\shortintertext{ \textbf{ourSolDep[5]} }
% \newBkt{\yb{1}}{\yb{2}}=    \yb{3},
% \newBkt{\yb{1}}{\yb{3}}=    \yb{4},
\newBkt{\yb{1}}{\yb{4}}&=(-\CC{1,4,2} \CC{2,3,4}) \yb{1}+\CC{1,4,2} \yb{2}
+\CC{1,4,3} \yb{3} & \qquad 
\newBkt{\yb{2}}{\yb{3}}&=\CC{2,3,4}   \yb{4} \\
\newBkt{\yb{2}}{\yb{4}}&=(-\CC{1,4,2}\CC{2,3,4}^2)\yb{1}+(\CC{1,4,2}\CC{2,3,
4})\yb{2}+(\CC{1,4,3}\CC{2,3,4})\yb{3} & \qquad 
\newBkt{\yb{3}}{\yb{4}}&=0 
\end{alignat*}
%% mySol[6] %%%%%%%%%%%%%%%%%%%%%%%%%%%%%%%%%%%%%%%%%%%%%%
%\textbf{ourSolDep[6]}
\begin{alignat*}{2}
\shortintertext{ \textbf{ourSolDep[6]} }
% \newBkt{\yb{1}}{\yb{2}} & =     \yb{3} \\
% \newBkt{\yb{1}}{\yb{3}} & =     \yb{4} \\
\newBkt{\yb{1}}{\yb{4}} & =\CC{1,4,3} \yb{3} & \qquad 
\newBkt{\yb{2}}{\yb{3}}& =\CC{2,3,1} \yb{1}+\CC{3,4,4} \yb{2}+\CC{2,3,4}
\yb{4}\\
\newBkt{\yb{2}}{\yb{4}}& =(\CC{1,4,3} \CC{2,3,4}+\CC{3,4,4}) \yb{3} & \qquad 
\newBkt{\yb{3}}{\yb{4}}& =(-\CC{1,4,3} \CC{2,3,1})\yb{1}+(-\CC{1,4,3} \CC{3,4,4}) \yb{2}+\CC{3,4,4} \yb{4}
\end{alignat*}

\section{About main maple scripts in this note}

Only Engel-try-1.mpl had a 2 pages full explanation,  
but 
Engel-mySols-v1.mpl,  DualAlg-Engel-v44.mpl, and  ZeroPlus-Engel-v2.mpl had
no chance to be explaned at all. 

In stead,  we invite readers to !http://math.akita-u.ac.jp/~mikami/!. 
There we prepare a tar-ball which includes revised maple scripts with 
the extension ".mpl.with", a small repository !big.mla!,  
and !please-read-me.txt!.

%\end{spacing}
\end{document}